\documentclass{article}
\usepackage[latin9]{inputenc}
\usepackage{multirow}
\RequirePackage{amsthm,amsmath,amsfonts,amssymb}
\RequirePackage{scalerel}
\RequirePackage[authoryear,round]{natbib} 
\RequirePackage{hyperref}
\hypersetup{
    colorlinks=true,
    citecolor=blue,
    linkcolor=blue}
\RequirePackage{graphicx}

\theoremstyle{plain}
\newtheorem{axiom}{Axiom}
\newtheorem{claim}[axiom]{Claim}
\newtheorem{theorem}{Theorem}[section]
\newtheorem{lemma}[theorem]{Lemma}
\newtheorem{proposition}[theorem]{Proposition}
\newtheorem{corollary}[theorem]{Corollary}

\theoremstyle{remark}

\def\bs{\boldsymbol}
\def\E{{\mathbb E}}
\def\P{{\mathbb P}}
\def\V{\mathrm{Var}}
\def\be{\begin{equation}}
\def\ee{\end{equation}}
\def\bea{\begin{eqnarray*}}
\def\eea{\end{eqnarray*}}
\def\bean{\begin{eqnarray}}
\def\eean{\end{eqnarray}}
\def\nn{\nonumber}
\def\nin{\noindent}

\def\ra{\rightarrow}
\def\Bl{\Bigl}
\def\Br{\Bigr}

\def\N{{\mathbb{N}}}

\def\X{{\bs{X}}}

\def\alp{\alpha}

\def\eps{\epsilon}

\def\Iapp{\mathcal{I}_{\mathrm{app}}}
\def\Capp{\mathcal{C}_{\mathrm{app}}}
\def\lmax{\ell_{\mathrm{max}}}
\DeclareMathOperator*{\bigtimes}{\scalerel*{\times}{\sum}}

\begin{document}

\title{Large-scale inference with block structure}
\author{Jiyao Kou and Guenther Walther\thanks{Email addresses: koujiyao@gmail.com, 
gwalther@stanford.edu. Supported by NSF grants DMS-1220311 and DMS-1501767}\\
        Stanford University}
\date{December 2021}

\maketitle

\begin{abstract}
The detection of weak and rare effects in large amounts of data arises
in a number of modern data analysis problems. Known results show that
in this situation the potential of statistical inference is severely
limited by the large-scale multiple testing that is inherent in these
problems. Here we show that fundamentally more powerful statistical
inference is possible when there is some structure in the signal that
can be exploited, e.g. if the signal is clustered in many small blocks,
as is the case in some relevant applications. We derive the detection boundary
in such a situation where we allow both the number of blocks and the
block length to grow polynomially with sample size. We derive these results
both for the univariate and the multivariate settings as well as for the problem of
detecting clusters in a network. These results recover
as special cases the sparse signal detection problem \citet{donoho2004higher}
where there is no structure in the signal, as well as the scan problem
\cite{chan2011detection} where the signal comprises a single interval.
We develop methodology that allows optimal adaptive detection in the
general setting, thus exploiting the structure if it is present without
incurring a relevant penalty in the case where there is no structure. The advantage
of this methodology can be considerable, as in the case of no structure the
means need to increase at the rate $\sqrt{\log n}$ to ensure detection,
while the presence of structure allows detection even if the means \emph{decrease} at a polynomial
rate.
\end{abstract}

\bigskip



\section{Introduction}

The problem of detecting a signal, such as an elevated mean, in a high-dimensional vector
of Gaussian observations has been of considerable interest
as it serves as the statistical model for the multiple
testing of a large number of hypotheses. 
This problem has been studied
in detail for the important setting where the signal is sparse and
weak, see Section~\ref{subsecHM}.
An important result of this research is that detection of the signal
is impossible unless the signal mean is at least of the order $\sqrt{\log n}$,
where $n$ is the sample size. This is a somewhat discouraging result in the context
of typical statistical inference problems, where a larger sample size usually allows to detect
a smaller mean. In fact, there is an earlier body of research that considers the above
detection problem in the case where the signal is aligned consecutively in an interval rather than
scattered at random. It can be shown that in
this ``block signal detection problem'' it is possible to detect much smaller means:
scan statistics can detect means that are sparse and weak and yet may \emph{decrease} at a rate
that is polynomial in $n$, see Section~\ref{sec:Comparison-with-other-statistic}
below. The stark contrast between these two results suggests that it may be
possible to perform statistical inference in the sparse and weak setting that is more powerful
in a fundamental and relevant way, provided there is some kind of structure in the signal
that can be exploited.

This paper develops methodology that is adaptive to such structure, i.e. it
automatically exploits structure that may be present in the data. We consider
a model where the signal is comprised of potentially many small blocks that are
scattered at random in the sequence, the \textquotedblleft
multiple blocks detection problem\textquotedblright . Two examples of such data are:

\begin{itemize}
\item[(i)] (Epidemic) Each location represents a one kilometer by one kilometer
square and the proportion of citizens that is diseased is measured
in each location. When there is a disease outbreak in the city, many
independent ``areas'' will have unusual high values, where ``area''
is defined as a two-dimensional block of locations. The task is to
detect whether there is a disease outbreak or not; see, for example,
\cite{kulldorff1999spatial,gangnon2001weighted,chan2009detection,walther2010optimal}.
In this example, the structure is ``spatial''.
\item[(ii)] (Financial) On each timestamp, we measure the predictive power
of a particular technical indicator for S\&P 500. Many periods with
unusual high predictive power indicate the potential usefulness of
the technical indicator for future trading, where ``period'' is
defined as a block of timestamps. The task is to detect whether the
technical indicator is useful or not. In this example, the structure
is ``temporal''.
\end{itemize}

In this paper we analyze the general setting of the multiple blocks
detection problem where both the number of blocks and the
block length can grow polynomially with the sample size. Note that
this model contains both the sparse signal detection problem
as well as the block signal detection problem as special cases. We establish
the detection boundary in this setting and introduce methodology that
allows optimal adaptive detection. That is, the methodology introduced
below will automatically utilize such structure if it is present and
provide optimal detection both when structure is present and when
it is not. Therefore this methodology is preferable whether prior information
about the number of blocks or the block length is available or not.

\subsection{Review of sparse signal detection\label{subsecHM}}

Consider an $n$-dimensional Gaussian vector with components
\be  \label{HM}
X_{i}=\mu\boldsymbol{1}_{\{\ell_1,\ldots,\ell_m\}}(i)+Z_{i},\qquad i=1,\ldots,n,
\ee
where the $Z_{i}$ are i.i.d. standard normal random variables and
the $m=n^{1-\beta}$, $0<\beta\le1$,
signal locations $\ell_1,\ldots,\ell_m$ are randomly
drawn from $\{1,2,\ldots n\}$ without replacement. 
We are testing whether $H_{0}:\mu=0$ vs $H_{1}:\mu=\mu(n)>0$. 

This \textit{sparse signal model} (or the closely related \textit{sparse heterogenous
mixture model} where $m \sim {\rm Bin}(n,n^{-\beta})$) has been investigated
by \cite{Ingster97,ingster1998minimax,IngsterSuslina,donoho2004higher,tony2011optimal}.
The extension to the case of dependent observations was investigated 
in \cite{delaigle2009higher,hall2010innovated,zhong2013tests}.
In a related context,
\cite{IngsterPouetTsy} consider the problem of  classifying a high-dimensional vector
from (\ref{HM}) as null or alternative, based on a training sample of several i.i.d. 
vectors from the alternative, and they derive a sharp classification boundary and classifiers
attaining this boundary.
\cite{VerzelenAC} consider the setting where one observes several i.i.d. copies of a high-dimensional
vector from a mixture of two Gaussians, assuming that the difference vector of the two Gaussian means 
is sparse. They consider
the problems of testing whether the difference in mean is zero and estimating which coordinates of
the difference are non-zero, and they derive minimax lower bounds and propose a number of methods
which attain these bounds.

In the testing problem for model (\ref{HM}) it turns out that there is a threshold effect for the 
likelihood ratio test. In the sparse regime where $\frac{1}{2}<\beta\le1$,
one calibrates $\mu=\mu(n)=\sqrt{2r\log n}$, where $0<r\le1$. By
\cite{Ingster97,ingster1998minimax,donoho2004higher}, the\emph{ detection boundary
}is defined as: 
\be  \label{HMboundary}
\rho^{*}(\beta)=\begin{cases}
\beta-\frac{1}{2} & \frac{1}{2}<\beta\le\frac{3}{4}\\
(1-\sqrt{1-\beta})^{2} & \frac{3}{4}<\beta\le1
\end{cases}
\ee
If $r>\rho^{*}(\beta)$, then $H_{0}$ and $H_{1}$ separate asymptotically,
i.e. the elevated mean can be detected with asymptotic probability one, while
if $r<\rho^{*}(\beta)$, then $H_{0}$ and $H_{1}$ merge asymptotically, i.e.
it is impossible to detect the elevated mean with power larger than the significance level.
Unfortunately, the likelihood ratio test requires a precise specification
of $r$ and $\beta$, so one would like to have a method which is adaptive
to the unknown $r$ and $\beta$ and perform as well as the likelihood
ratio test.

\cite{IngsterSuslina} introduced such an adaptive method by combining three
different procedures, while \cite{donoho2004higher} proposed the following \textit{higher
criticism} (HC):

\be  \label{HC}
HC_{n}=\max_{1\le i\le\frac{n}{2}}\sqrt{n}\frac{\frac{i}{n}-p_{(i)}}{\sqrt{p_{(i)}(1-p_{(i)})}},
\ee
where $p_i :=\bar{\Phi} (X_i)$ is the p-value for $X_i$ and the $p_{(i)}$ denote
the p-values sorted in increasing order.
It can be shown that HC can attain the detection boundary so it is
optimal for sparse signal detection: HC will separate
the two hypotheses asymptotically whenever the likelihood ratio test
can asymptotically separate the two hypotheses.

Another popular choice is called \textit{Berk-Jones} statistic (BJ),
which is defined as follows:

\be  \label{BJ}
BJ_{n}=\max_{1\le i\le \frac{n}{2}}\left(i\log\frac{i}{np_{(i)}}+
(n-i)\log\frac{1-\frac{i}{n}}{1-p_{(i)}}\right)\mathbf{1}
\left(p_{(i)}<\frac{i}{n}\right),
\ee
BJ is also optimal for sparse signal detection,
see \cite{donoho2004higher}, and its finite sample performance appears
to be better than that of HC, see \cite{walther2013average,li2015higher}.

Some alternative tests were studied in \cite{jager2007goodness,zhong2013tests,walther2013average}.
In particular, \cite{jager2007goodness} have shown that all members of the
$\phi$-divergence family
$S_{n}^{+}(s), s\in[-1,2]$, attain the detection boundary (\ref{HMboundary}),
where $ S_n^{+}(s)=n\,\max_{1\le i\le\frac{n}{2}}
K_{s}\left(\frac{i}{n},p_{(i)}\right) \mathbf{1} \left(p_{(i)}<\frac{i}{n}\right)$ 
and $K_s(\cdot,\cdot)$ is given in \cite{jager2007goodness}.
This family  contains as special cases the Berk-Jones
statistic ($s=1$) and for $s=2$ a statistic that is equivalent to the higher criticism:
$S_n^+(2)=\frac{1}{2} (HC_n^+)^2$.

\cite{IngsterSuslina,tony2011optimal} extend the detection boundary to the case
$0<\beta<\frac{1}{2}$ which  \cite{tony2011optimal} call the \textit{ dense regime} (the
designation \textit{ moderately sparse} in \cite{IngsterPouetTsy} is perhaps more apt), 
and they show that there is also a threshold
effect for the likelihood ratio test. In the dense case one needs to calibrate
$\mu=\mu(n)=n^{r}$. Then the detection boundary is defined as:
\be  \label{HMdenseBoundary}
\rho^{*}(\beta)=\beta-\frac{1}{2}
\ee
If $r>\rho^{*}(\beta)$, $H_{0}$ and $H_{1}$ separate asymptotically
and if $r<\rho^{*}(\beta)$, $H_{0}$ and $H_{1}$ merge asymptotically.
It is shown in \cite{tony2011optimal} that HC is also optimal for
detection in the dense case. Note that the dense
case is much less challenging since even a simple z-test will do very well,
see \cite{IngsterSuslina}.

\subsection{Organization of the paper and notation}

In Section~\ref{sec:Def-Block-heterogeneous-mixtures}
we introduce the definition of the multiple blocks model
and derive the detection boundary for this model.
In Section~\ref{sec:Spatial-HC-BJ-null} we propose procedures for detection
in this model, namely the
structured higher criticism and structured Berk-Jones statistics, and more
generally the family of structured $\phi$-divergences, and we
evaluate their properties under the null distribution. This section also derives
tail bounds for the higher criticism and Berk-Jones statistics which may
be of independent interest. In Section~\ref{sec:Optimality-of-spatial-HC-BJ},
we establish the optimality of these  statistics for the multiple blocks model.
In Section~\ref{sec:Comparison-with-other-statistic}, we compare the
performance of structured higher criticism and structured Berk-Jones
statistics with other methods. In Section~\ref{sec:Simulation}, a
simulation study is carried out to illustrate our results.
Section~\ref{multivariate} treats the multivariate case, and Section~\ref{graphs}
deals with clusters in a network.
Section~\ref{composite} addresses composite alternatives.
In Section~\ref{sec:Conclusion-and-Discuss}, we discuss some
possible extensions and future research topics. All proofs of the main
theorems and propositions are put in Section~\ref{sec:Proofs}.
Some technical arguments are deferred to the \hyperref[appn]{Appendix}.

We denote the number of design points contained in a set $I$ by $|I|$.
For the half-open intervals and rectangles we consider here this will
typically be equal to the Lebesgue measure of $I$.
$L_{n}$ denotes terms satisfying $\log L_{n}=o(\log n)$,
which may vary from place to place. Note that for all fixed $\epsilon>0$,
$L_{n}n^{\epsilon}\rightarrow\infty$ and $L_{n}n^{-\epsilon}\rightarrow0$
as $n\rightarrow\infty$. We employ the usual $O_p$ and $o_p$ notation for
a sequence of random variables $X_n$ and in
addition write $X_{n}=\Omega_{p}(a_{n})$ if for every $\epsilon \in (0,1)$ there
exists a finite $M>0$ such that $P(|X_{n}/a_{n}|<M)<\epsilon$ for
all $n$ that are large enough. 
In this paper, $\log n$ is used for the natural logarithm
while $\log_{2}n$ is used for logarithm to the base 2.

\section{The multiple blocks model\label{sec:Def-Block-heterogeneous-mixtures}}

The sparse signal model (\ref{HM}) posits
that there is no structure in the signal. However, it turns out that if some
structure does exist, then the detection problem becomes easier in a fundamental way
and a much better result is attainable. Specifically, in this paper, we consider
the situation where the signal is clustered into multiple blocks  with unknown
length. We call this the multiple blocks  model:

\begin{equation}
X_{i}=\mu\boldsymbol{1}_{\bigcup_{g=1}^{m}I_{g}}(i)+Z_{i},\qquad i=1,\ldots,n
\label{eq: Def-Block-mixture-problem}
\end{equation}
where the $I_g$ are mutually disjoint intervals at random locations
 and the $Z_{i}$ are i.i.d. standard normal random variables. 
The difficulty of this detection problem depends on the size of $\mu$,
the number $m$ of blocks, and the minimum block length $\min_g |I_g|$.
In order to derive a succinct theoretical result about the detection boundary
we let the number of blocks $m=n^{1-\alpha-\beta}$ and assume that each of the unknown
blocks $I_{g}$ has equal length $|I_{g}|=n^{\alpha}$, $g=1,\ldots,m$,
where $0\le\alpha<1$ and $0<\alpha+\beta\le1$. 
The task is to test $H_{0}:\mu=0$
vs. $H_{1}:\mu=\mu(n)>0$. All of the following results can be reformulated
for unequal block lengths in terms of $\min_g |I_g|$ by using a minimax
statements for lower bounds.

If $\alpha=0$, then $|I_{g}|=1$ for all $g=1,\ldots,m$ and we obtain
the sparse signal model (\ref{HM}). If $\alpha=1-\beta$, then
we only have $m=1$ block, and our problem reduces to the block
signal detection model (\ref{eq: def-block-signal-model}) discussed below. 
Thus our multiple blocks model
is a generalization of both the sparse signal detection problem and
the block signal detection problem.

Another special case of the multiple blocks model is investigated by
\cite{jeng2010optimal}. They use a likelihood ratio selection procedure
for detecting very sparse and very short
segments of elevated means, i.e both the number and the lengths of
segments grow at most logarithmically with sample size.
Their results suggest that this likelihood procedure will not be able
to attain the detection boundary for the more general model considered here.

The multiple blocks model describes a situation where
the signal arises in many locations in the form of small clusters.
While this model can be analyzed with HC or BJ, the results in 
Sections~\ref{subsecHM} and \ref{sec:Comparison-with-other-statistic} suggest that such an
analysis would be quite suboptimal: In the sparse case $\beta > \frac{1}{2}$,
HC and BJ require that each of
the $n^{1-\beta}$ signal means is of size at least $\sqrt{c(\beta) \log n}$
for some constant $c(\beta)$. In contrast, if the $n^{1-\beta}$
signal means are aligned in one single interval, then a certain scan statistic
will detect signal means as small as $\frac{\sqrt{2 \beta \log n}}{
n^{(1-\beta)/2}}$, which is a drastically smaller threshold, see 
Section~\ref{sec:Comparison-with-other-statistic}.
 This is due to the square root
law which the scan exploits in this situation. These results
suggest that likewise in the multiple blocks model it
should be possible to drastically improve upon the power of HC and BJ
by exploiting the structure of the signal. It will be shown below
how this can be done by introducing the \textit{structured} HC and BJ statistics.
To this end, we first derive the detection boundary for this problem.

\subsection{The detection boundary for the multiple blocks model}

As in the sparse signal detection problem, the calibration of the detection
boundary differs in the sparse and in the dense case, but the sparse
case is now defined by the condition 
$\frac{\beta}{1-\alpha}>\frac{1}{2}$.
\begin{theorem} \label{thm: Block-mixture-detection-boundary}
Consider the multiple blocks model \eqref{eq: Def-Block-mixture-problem}.
\begin{itemize}
\item [(i)](Sparse case) If
$\frac{\beta}{1-\alpha}>\frac{1}{2}$ set
\begin{equation}
\mu=\mu(n)=\sqrt{2r\log n}/\sqrt{n^{\alp}}\label{eq: sparse-calibrate}
\end{equation}
with $r>0$ and
\begin{equation}
\rho^{*}(\alpha,\beta)=\begin{cases}
\beta-(1-\alpha)/2 & \text{ if }\beta/(1-\alpha)<\frac{3}{4}\\
(\sqrt{1-\alpha}-\sqrt{1-\alpha-\beta})^{2} & \text{ if }\beta/(1-\alpha)\ge\frac{3}{4}
\end{cases}\label{eq: detection_boundary_sparse}
\end{equation}
\item [(ii)](Dense case) If $\frac{\beta}{1-\alpha}<\frac{1}{2}$ set
\begin{equation}
\mu=\mu(n)=n^{r}/\sqrt{n^{\alp}}\label{eq: dense-calibrate}
\end{equation}
and
\begin{equation}
\rho^{*}(\alpha,\beta)=\beta-\frac{1-\alpha}{2}.\label{eq: detection_boundary_dense}
\end{equation}
\end{itemize}

If $r<\rho^{*}(\alpha,\beta)$, 
then $H_{0}$ and $H_{1}$ merge asymptotically, i.e. the sum of
Type I and Type II errors tends to 1 for any test.
\end{theorem}

The proof of the theorem can be derived from (\ref{HMboundary}) and (\ref{HMdenseBoundary})
by applying the square root law.
The theorem shows that $\rho^{*}(\alpha,\beta)$
is a lower bound for model \eqref{eq: Def-Block-mixture-problem}:
If $r<\rho^{*}(\alpha,\beta)$, then detection
is not possible. In Sections~\ref{sec:Spatial-HC-BJ-null} and 
\ref{sec:Optimality-of-spatial-HC-BJ} we will derive and investigate
procedures which attain this lower bound when both the sparsity
level and the block length are unknown, i.e. these procedures are adaptive to
both $\alpha$ and $\beta$. Hence $\rho^{*}(\alpha,\beta)$
does in fact describe the detection boundary for the 
model \eqref{eq: Def-Block-mixture-problem}.

Note that the calibration of $\mu$ in Theorem~\ref{thm: Block-mixture-detection-boundary}
has the divisor $\sqrt{n^{\alp}}$ which does not appear in the sparse signal
detection problem. This shows that in the sparse case the multiple blocks model allows 
the detection of much smaller means. Even if the blocks are very
short, say of length 2 or 3, this will improve upon the detection boundary
(\ref{HMboundary}). Longer blocks, e.g. of order $\log n$ or $n^{\alpha}$,
have an even more dramatic effect by changing the scaling of the
detection boundary. It is interesting to note that in the dense case 
$\mu^{*}:=n^{\rho^{*}(\alpha,\beta)}/\sqrt{n^{\alp}}=n^{\beta-\frac{1}{2}}$
does not depend on $\alpha$, which suggests that
the block structure may not be important anymore in the dense case.
We discuss this issue further in Section~\ref{subsec:Discussion:Sparse-To-Dense}.

\section{The structured higher criticism and Berk-Jones 
statistics\label{sec:Spatial-HC-BJ-null}}

In order to motivate our approach we note that detection in the multiple
blocks model requires to aggregate the evidence in the data in two
ways: For a given candidate interval the evidence must be combined
within that interval, e.g. by averaging the observations. Then
this evidence must be aggregated across intervals by a multiple
testing procedure such as HC. However, a straightforward implementation of
this idea is not promising: The detection
boundary \eqref{HMboundary} in the unstructured case is due to the 
multiple testing of $n$ p-values. If one were to compute a p-value
for each candidate interval, then the  ensuing
massive multiple testing problem results in about $ n^2$ p-values
and HC may not attain the detection boundary 
\eqref{eq: detection_boundary_sparse}.
Moreover, many of these p-values will be highly correlated and so
the usual critical values for HC are not applicable.

We circumvent these problems by considering
an appropriate approximating set of intervals that possesses
the following three properties: First, each of the about $ n^2/2$
intervals with endpoints in $\{1,\ldots,n\}$ can be approximated
sufficiently well by an interval in the approximating set so that
the resulting approximation error to the signal does not detract from
the detection boundary.
 Second, there are only $O(n \log n)$ intervals
in the approximating set. As a consequence, the multiple testing
does not become noticeably more difficult as HC still has to 
assess only of the order $n$ p-values rather than $n^2$. Third, the
approximating set is sparse enough to allow an analysis of the null
distribution of HC in the context of independent p-values, as will be explained below.

These criteria are satisfied by the approximating set used in
 \cite{walther2010optimal,rivera2013optimal}: 
\smallskip

For each level $\ell=0,\ldots,\lmax$, where 
$\lmax=\text{\ensuremath{\lceil}}\log_{2}\frac{n}{8}\rceil$:

\[
\Iapp (\ell):=\Bigl\{ (j,k] \subset (0,n]:
j,k\in\{id_{\ell},i=0,1,\ldots\}\text{ and }2^{\ell-1}<k-j\le2^{\ell}\Bigr\} 
\]
where $d_{\ell}=\lceil \eps_{\ell} 2^{\ell-1} \rceil$ for 
$\eps_{\ell}=\frac{1}{6\sqrt{\log_{2}\frac{n}{2^{\ell-1}}}}=\frac{1}{6\sqrt{\lmax-\ell+4}}$.

That is, the collection $\Iapp (\ell)$ approximates intervals with lengths in $(2^{\ell -1},2^{\ell}]$
via endpoints on a grid whose spacing is a fraction $\eps_{\ell}$ of the approximate interval
length $2^{\ell -1}$, where the precision parameter $\eps_{\ell}$ changes with the length of the
intervals such that it produces a finer approximation for smaller intervals.
The approximating set $\bigcup_{\ell} \Iapp (\ell)$ has cardinality $O(n \log n)$ but
approximates all intervals sufficiently well to allow optimal inference, see 
Proposition~\ref{properties} in Section~\ref{sec:Proofs}
for a more precise statement of its properties.

Now we define \textit{structured higher criticism} $sHC_{n}$
and \textit{structured Berk-Jones} statistic $sBJ_{n}$ as follows:
\[
sHC_{n}\ =\ \max_{\ell=0}^{\lmax}\sqrt{\frac{n}{2^{\ell}n_{\ell}}}HC_{n_{\ell}}(\ell),
\]

\[
sBJ_{n}\ =\ \max_{\ell=0}^{\lmax}\,\frac{n}{2^{\ell}n_{\ell}}BJ{}_{n_{\ell}}(\ell),
\]
where $HC_{n_{\ell}}(\ell)$ and $BJ{}_{n_{\ell}}(\ell)$ denote the one-sided higher 
criticism (\ref{HC}) and Berk-Jones statistic (\ref{BJ})
evaluated on the 
$$
n_{\ell}\ :=\ \#\Iapp(\ell)
$$ 
p-values
pertaining to $\Iapp(\ell)$, i.e. the p-values $\{\bar{\Phi} (\X(I)),\, I\in \Iapp(\ell)\}$,
where $\X(I):=\sum_{i\in I}X_{i}/\sqrt{|I|}$ is the standardized average
over the interval $I$.

More generally, we define for $s \in [-1,2]$ the \textit{structured $\phi$-divergence}
\[
sS_n^{}(s)\ =\ \max_{\ell=0}^{\lmax}\,\frac{n}{2^{\ell}n_{\ell}}S_{n_{\ell}}^{+}(s,\ell)
\]
where likewise $S_{n_{\ell}}^{+}(s,\ell)$ denotes the $\phi$-divergence $S_{n_{\ell}}^{+}(s)$ defined 
in Section~\ref{subsecHM}, evaluated on the p-values that
pertain to $\Iapp (\ell)$. 
In particular, $sS_n^{}(1)=sBJ_{n}$
and $sS_n^{}(2)$ is equivalent to $sHC_{n}$.

The difficulty in analyzing the null distributions of $BJ{}_{n_{\ell}}(\ell)$
and $HC_{n_{\ell}}(\ell)$ lies in the fact that the underlying p-values
are no longer independent because they are based on data pertaining
to intervals that may overlap. The key to controlling 
those null distributions is the sparse construction of the approximating set
$\Iapp (\ell)$:
It is shown in Lemma~\ref{lem: Short-range-dependent} that
the intervals in $\Iapp (\ell)$
can be grouped into a small number of groups such that each group
contains about $\frac{n}{2^{\ell}}$ intervals that are disjoint and whose
corresponding p-values are therefore independent.
Hence the empirical measure of the p-values can
be written as an average of a small number of empirical measures,
each of which is based on independent p-values. This allows to
use Jensen's inequality to bound $BJ{}_{n_{\ell}}(\ell)$ and $HC_{n_{\ell}}(\ell)$
by the maximum of a small number of such statistics, each
of which is based on $\frac{n}{2^{\ell}}$ independent p-values. 
This maximum can then be controlled
via tail bounds for these statistics. 
Furthermore, this explanation shows that the scaling in
$HC_{n_{\ell}}(\ell)$ should be $\sqrt{\frac{n}{2^{\ell}}}$
rather than $\sqrt{n_{\ell}}$, 
hence the rescaling factor $\sqrt{\frac{n}{2^{\ell}n_{\ell}}}$ 
for $HC_{n_{\ell}}(\ell)$, and analogously for $BJ{}_{n_{\ell}}$ and the structured
$\phi$-divergence.

\begin{theorem}
\label{thm:HC-BJ-Under-the-null}Under the null hypothesis $\mu=0$,
\begin{align*}
\frac{sBJ_{n}}{\log\log n}\  & \stackrel{p}{\leq}\ 3\ (n\rightarrow\infty)\\
sHC_{n}\  & =\ O_{p}(\log^{2}n)
\end{align*}
\end{theorem}

Note that under the null distribution $\frac{BJ_{n}}{\log\log n}\stackrel{p}{\rightarrow}1$
and $\frac{HC_{n}}{\sqrt{2\log\log n}}\stackrel{p}{\rightarrow}1$,
see \cite{jager2007goodness}. Thus the penalty for additionally examining
structure in the data is at most a factor of 3 for $sBJ_{n}$.
In particular, the more general $sBJ_{n}$ is still optimal
in the special case (\ref{HM}) when there is no structure
in the signal, and likewise for $sHC_{n}$.

As an aside, it is not clear that the result of Theorem~\ref{thm:HC-BJ-Under-the-null}
for $sHC_{n}$ can be improved
as the smallest p-values have heavy tails, see \cite{walther2013average}.
While that can be controlled in the case of a single HC statistic, see
e.g. page 601-603 of \cite{shorackW}, $sHC_{n}$
is the maximum of $\sim\log n$ terms that involve HC statistics. 
In the context of a sparse Gaussian graphical model
\cite{FanJinYao} give a result about how to group correlated p-values to guarantee independence
within each group, but this result is not applicable here.

For the proof of the theorem we will need the following tail bounds
which may be of independent interest: 
\begin{proposition}
\label{prop:HC-BJ-each-lvel-null}
Let $F_{n}$ be the empirical cdf
of $U_{1},\ldots,U_{n}$ i.i.d. $U(0,1)$ and let $U(\cdot)$ be a
standard Brownian Bridge. For $0<a<b<1$ and $\eta>0$: 

\begin{itemize}
\item[(i)]
\[
\P\left(\sup_{t\in[a,b]}\frac{U(t)}{\sqrt{t(1-t)}}>\eta\right)\ \le \ 
\frac{\frac{2}{\eta}+\eta\log\frac{b(1-a)}{a(1-b)}}{\sqrt{2\pi}}\ e^{-\eta^2/2}
\]

\item[(ii)]
\begin{align*}
  & \P\left(\sup_{t\in[a,b]}n\left(F_{n}(t)\frac{\log F_{n}(t)}{t}+(1-F_{n}(t))
  \log\frac{1-F_{n}(t)}{1-t}\right)>\eta\right)\\
 &  \le \  2e\left(\eta\log \frac{b(1-a)}{a(1-b)}+1\right)\exp(-\eta)
\end{align*}

\item[(iii)] For every $K>1$: 
\[
\P_{H_{0}}(BJ_{n}>\eta)\ \le \ 22K(\log n)(\eta+1)\exp(-\eta)+2n^{1-K}
\]

\item[(iv)]
\[
\P_{H_{0}}(HC_{n}>\eta)\ \le \ \P\left(\sup_{t\in(0,1)}\sqrt{n}\frac{F_{n}(t)-t}{
\sqrt{t(1-t)}}>\eta\right)\ \le\ \frac{C}{\eta}
\]
\hspace{0.5in} for $\eta\ge\sqrt{D\log\log n}$ with $D>2$ where the constant C
depends only on D.
\end{itemize}
\end{proposition}
\nin \cite{miller1982maximally} give a two-sided bound corresponding
to (i) which holds asymptotically.
(ii) improves the exponential bound provided in
 \cite{duembgen2014confidence}.
As for (iv), there exists no exponential
inequality for $HC_{n}$ due to the heavy algebraic tails of the smallest
p-values, see \cite{walther2013average}.

\section{Optimality of the structured higher criticism and structured Berk-Jones
statistic for the multiple blocks model\label{sec:Optimality-of-spatial-HC-BJ}}

The following theorem shows that every structured $\phi$-divergence,
and in particular the structured higher criticism
and the structured Berk-Jones statistic, attain the lower bound 
established in Theorem \ref{thm: Block-mixture-detection-boundary}
for the sparse case. The theorem also shows that the structured higher
criticism statistic attains the lower bound in the dense case.
Thus these procedures are
optimal for detection in the multiple blocks model and are adaptive to both
the unknown block length and the unknown sparsity level.

\begin{theorem}
\label{thm: Spatial HC is optimal}
Consider the multiple blocks model \eqref{eq: Def-Block-mixture-problem}.
\begin{itemize}
\item[(i)] In the sparse case $\frac{\beta}{1-\alpha}>\frac{1}{2}$ with the calibration
\eqref{eq: sparse-calibrate} for the mean of the signal, let
$r>\rho^{*}(\alpha,\beta)$ in \eqref{eq: detection_boundary_sparse}.
Then every member of the family of structured $\phi$-divergences
$sS_n^{}(s)$, $s \in [-1,2]$, has asymptotic power 1. 
\item[(ii)] In the dense case $\frac{\beta}{1-\alpha}<\frac{1}{2}$ with the calibration
\eqref{eq: dense-calibrate} for the mean of the signal, let
$r>\rho^{*}(\alpha,\beta)$ in \eqref{eq: detection_boundary_dense}.
Then $sHC_{n}$ has asymptotic power 1.
\end{itemize}
\end{theorem}

We note
that while there are $O(n^{2})$ possible intervals, the use
of the approximation set makes it possible to compute these structured statistics
in $O(n\log^{2}n)$ time, almost linear in the
number of observations.

As a corollary to the above theorem we note that $sHC_{n}$ and $sBJ_{n}$ are
optimal for sparse signal  detection and for block signal
detection, which are special cases the model \eqref{eq: Def-Block-mixture-problem}:
\begin{corollary}
$sHC_{n}$ and $sBJ_{n}$ achieve the optimal detection
boundary (\ref{HMboundary}) in the sparse signal model
(\ref{HM}).
\end{corollary}

The corollary follows upon observing that the sparse signal model (\ref{HM})
obtains as the
special case $\alpha=0$. By Theorem \ref{thm: Spatial HC is optimal}(i),
$sHC_{n}$ and $sBJ_{n}$ can reliably detect the alternative
if $r>\rho^{*}(0,\beta)$ in \eqref{eq: detection_boundary_sparse},
 which equals the detection boundary (\ref{HMboundary}). 

In the block signal detection problem the signal is aligned in an interval $I_n$, i.e.
\begin{equation}
X_{i}=\mu\mathbf{1}_{I_{n}}(i)+Z_{i},\qquad i=1,\ldots,n\label{eq: def-block-signal-model}
\end{equation}

\begin{corollary}
$sHC_{n}$ and $sBJ_{n}$ achieve the optimal detection
boundary for the block signal detection problem (\ref{eq: def-block-signal-model}).
\end{corollary}

The block signal detection problem corresponds to $\alpha=1-\beta$.
By Theorem \ref{thm: Spatial HC is optimal}(i), $sBJ_{n}$ and
$sHC_{n}$ can reliably detect the alternative if $r>\rho^{*}(1-\beta,\beta)$ in
\eqref{eq: detection_boundary_sparse}, where

\[
\rho^{*}(1-\beta,\beta)=\beta=1-\alpha.
\]
Thus, when writing the alternative in terms of $\mu$, we can reliably detect the
alternative if 
\[
\mu>(1+\epsilon)\sqrt{2(1-\alpha)\log n}/\sqrt{n^{\alp}}=
(1+\epsilon)\sqrt{2\log\frac{n}{|I_n|}}/\sqrt{|I_n|},
\]
for any $\epsilon>0$, which matches the optimal detection boundary
for block signal detection given in Section~\ref{sec:Comparison-with-other-statistic} in terms of rate
and constant. (The more refined result in Section~\ref{sec:Comparison-with-other-statistic} even allows
$\epsilon_n \downarrow 0$ at a certain rate for the penalized scan,
and it is not clear whether $sBJ_{n}$ or
$sHC_{n}$ can attain that behavior near the boundary.)

\section{Comparison with other
methods\label{sec:Comparison-with-other-statistic}}

In this section we compare structured BJ and HC with relevant other
methodology in terms of their theoretical performance. Section~\ref{sec:Simulation}
will complement this comparison with a simulation study.

Perhaps the most obvious approach to the multiple blocks model
is to directly use HC or BJ. Note that this approach ignores the block
structure in the data. 

In the sparse unstructured case $\frac{1}{2}<\beta\le1$, if we use the calibration 
\eqref{eq: sparse-calibrate}, then the detection boundary (\ref{HMboundary}) for HC
becomes
\[
\rho_{HC}^{*}(\alpha,\beta)=\begin{cases}
(\beta-\frac{1}{2}) n^{\alp} & \frac{1}{2}<\beta<\frac{3}{4}\\
(1-\sqrt{1-\beta})^{2} n^{\alp} & \beta\ge\frac{3}{4}.
\end{cases}.
\]
In the dense unstructured case $0<\beta<\frac{1}{2}$, if we use the calibration 
\eqref{eq: dense-calibrate}, then the detection boundary (\ref{HMdenseBoundary})
for HC becomes
\[
\rho_{HC}^{*}(\alpha,\beta)=\beta-\frac{1-\alpha}{2}.
\]

While the above detection boundaries are for the unstructured case,
it follows that HC and BJ cannot improve on these boundaries in the
multiple blocks model because they are invariant under permutations
of the observations and hence the block structure has no effect on
the inference.
Therefore HC and sHC compare as follows:

\begin{itemize}
\item[1.] When $\beta>\frac{1}{2}$ (and so $ $$\frac{\beta}{1-\alpha}>\frac{1}{2}$),
then both HC and sHC are in the sparse regime. Compared to sHC, the detection
boundary for HC is increased by a factor of $\sqrt{n^{\alp}}$.
Unless $\alpha=0$ (i.e. the length of the block is 1), the loss of
power of HC is significant.

\item[2.] When $\frac{1-\alpha}{2}<\beta<\frac{1}{2}$, then HC is in the dense
regime and sHC is in the sparse regime. Nevertheless, sHC has a more
favorable detection boundary: Compared to sHC, the detection
boundary for HC is increased by a factor (up to a $\log n$ term)
of $\frac{\sqrt{n^{\alp}}}{n^{\frac{1}{2}-\beta}}=n^{\beta-\frac{1-\alpha}{2}}$,
which grows polynomially with $n$. Therefore the loss of power of HC
is also significant.

\item[3.] When $\frac{}{}$$\beta<\frac{1-\alpha}{2}$ (and so $\beta<\frac{1}{2}$),
then both HC and sHC are in the dense regime. The detection boundaries
are the same for both methods and thus both HC and sHC are optimal
for the multiple blocks model. The reason for this is
that now the fraction of elevated means is so large that 
the block structure does not provide a noticeable benefit any more.
\end{itemize}

In light of the block structure in the data, another alternative approach
would be to use a scan statistic. Note that a scan statistic is designed
to detect a signal on an interval but not to aggregate the evidence
across multiple intervals. It is shown in \cite{chan2011detection}
that the scan with scale-dependent critical values, such as the
penalized scan 

\be  \label{penscan}
P_{n}=\max_{0\le j<k\le n}\left(\frac{\sum_{i=j+1}^{k}X_{i}}{\sqrt{k-j}}-\sqrt{2\log\frac{en}{k-j}}
\right)
\ee
dominates the regular scan, so we 
will only discuss the former. Moreover, it is shown in \cite{chan2011detection}
that evaluating the penalized scan on an approximating set:
\be   \label{apprpenscan}
P_{n}^{\mathrm{app}}=\max_{I\in\bigcup_{\ell=0}^{\lmax}
\Iapp(\ell)}\left(\frac{\sum_{i\in I}X_{i}}{\sqrt{|I|}}
-\sqrt{2\log\frac{en}{|I|}}\right)
\ee 
will not detract from its performance, while reducing the computational
effort from $O(n^2)$ to $O( n \log n)$. 

For the block signal detection problem (\ref{eq: def-block-signal-model}) 
where the signal is aligned in an interval $I_n$,
it is shown in
\cite{chan2011detection} that $P_{n}^{\mathrm{app}}$ has asymptotic power one if 
$\mu =\mu(n) \ge(\sqrt{2}+\epsilon_n)\sqrt{\log\frac{en}{|I_{n}|}}/\sqrt{|I_{n}|}$ 
with $\epsilon_n \sqrt{\log\frac{en}{|I_{n}|}} \rightarrow \infty$, while no
consistent test exists if $\mu=\mu(n)\le(\sqrt{2}-\epsilon_n)\sqrt{\log\frac{en}{|I_{n}|}}/
\sqrt{|I_{n}|}$ with $\epsilon_n \sqrt{\log\frac{en}{|I_{n}|}} \rightarrow \infty$.
Thus $P_{n}$ and $P_{n}^{\mathrm{app}}$ are optimal tests if the signal is aligned in a single interval. 
See \cite{walther2010optimal} for  a corresponding result
in the multivariate case and \cite{arias2005near} for earlier work
deriving the threshold $\sqrt{2\log n}/\sqrt{|I_{n}|}$ for the regular (unpenalized)
scan, which is optimal for very short interval lengths $|I_n|$ up to about $\log n$.

If we consider instead the multiple blocks model (\ref{eq: Def-Block-mixture-problem}),
then we obtain the following result:

\begin{theorem}
\label{thm:detection-boundary-penalized-scan}
In the multiple blocks model (\ref{eq: Def-Block-mixture-problem})
the detection boundary for the penalized scans $P_n$ and $P_{n}^{\mathrm{app}}$
given in (\ref{penscan}) and (\ref{apprpenscan}) is
$$
\rho_{\mathrm{pen}}^{*}(\alpha,\beta)=
\begin{cases}
\beta-(1-\alpha)/2 & \text{ if }\beta/(1-\alpha)<\frac{1}{2},\\
(\sqrt{1-\alpha}-\sqrt{1-\alpha-\beta})^{2} & \text{ if }\beta/(1-\alpha)>
\frac{1}{2},
\end{cases}
$$
with calibration \eqref{eq: dense-calibrate} in the first case and
calibration \eqref{eq: sparse-calibrate} in the second.
\end{theorem}

Thus the penalized scan attains the optimal detection boundary
except in the case $\frac{1-\alpha}{2}<\beta<\frac{3(1-\alpha)}{4}$,
where  $\rho_{\mathrm{pen}}^{*}(\alpha,\beta)$ is larger than
$\rho^{*}(\alpha,\beta)$ given in (\ref{eq: detection_boundary_sparse}).

\subsection{\label{subsec:Discussion:Sparse-To-Dense}Discussion: What matters
for good inference?}

Efficient inference in the multiple blocks model requires to combine
the evidence in two different ways: the evidence within a block needs
to be combined in order to make use of the square root law, and then
this evidence needs to be aggregated across blocks.

In the very sparse case $\frac{\beta}{1-\alpha}\ge\frac{3}{4}$,
the block structure is the most important aspect. In order to aggregate
 the information across blocks it is sufficient to simply scan for
the maximum of the within-block statistics. For this reason,
the penalized scan and sHC/sBJ perform well, whereas HC and BJ exhibit
a severe loss of power because they do not make use of the structure in
the signal and therefore forego the considerable advantage that derives
from the square root law.

In the  sparse case $\frac{1}{2}<\frac{\beta}{1-\alpha}<\frac{3}{4}$,
the block structure is still very important. However, optimally
aggregating the information across blocks requires an approach that
is more sophisticated than simply scanning for the maximum of the within-block
statistics. This explains why sHC/sBJ are optimal while
HC and BJ still exhibit a severe loss
of power as they do not make use of the structure in the signal.

The dense case $\frac{\beta}{1-\alpha}<\frac{1}{2}$ turns out to be the
regime where the structure in the signal is of no help for inference any
more. The reason for this perhaps surprising fact is that the fraction
of elevated means is now so large that asymptotic optimality obtains via the
square root law by simply averaging all observations, i.e. performing
a z-test. While HC and BJ are geared towards the sparse
case, they do attain the detection boundary in this dense case also, and
so do the structured versions sHC and sBJ and the penalized scan.

\section{Simulation study\label{sec:Simulation}}

This section  provides a simulation study that compares the
performance of sHC, sBJ, HC, BJ, and the penalized scan.
The sample size is $n=10000$ and power is with respect to a significance
level of $5\%$. Critical values for this significance level were simulated
with 10000 simulations and power was estimated with 2000 simulations.

\subsection{Simulation results for the very sparse case\label{subsec:Simulation-result-for-very-sparse}}

\begin{figure}
\includegraphics[height=0.43\textheight]{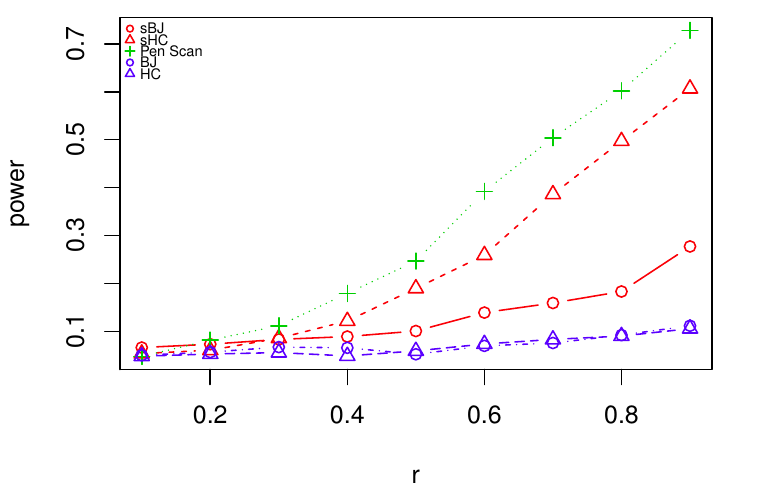}
\caption{Power for sBJ (solid line with $\circ$), sHC (short dashed line with $\triangle$), BJ (dash-dot line with
$\circ$), HC (long dashed line with $\triangle$) and penalized scan (dotted line with $+$) 
in the very sparse case $\alpha=0.2$, $\beta=0.65$.
\label{fig:Comparison-of-Powers-very-sparse}}
\end{figure}

We set $\alpha=0.2$
and $\beta=0.65$, so $\beta/(1-\alpha)\approx0.813$.
Power for the various methods is plotted in
 Figure~\ref{fig:Comparison-of-Powers-very-sparse}
as a function of $r$ in the calibration (\ref{eq: sparse-calibrate}).
The plot shows that the penalized scan 
has the highest power, followed by the structured HC and sBJ.
HC and BJ are nearly powerless even for large values of $r$. 
This simulation result confirms our conclusions
from Section~\ref{sec:Comparison-with-other-statistic}.
sBJ has less power than sHC partly because the first few p-values
in the appropriate level contain the most information in the very sparse
regime and sHC effectively puts more weights toward those  than
sBJ, see \cite{walther2013average} for an explanation of this phenomenon
in the setting without structure.

\subsection{Simulation results for the sparse case}

\begin{figure}
\includegraphics[height=0.43\textheight]{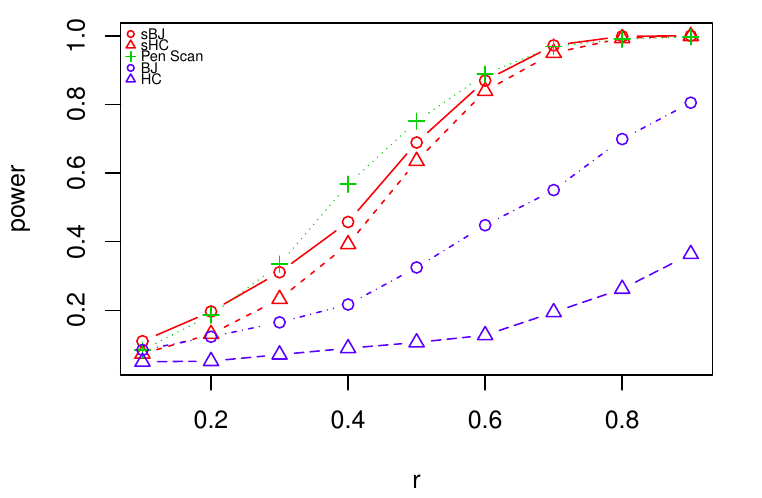}
\caption{Power for sBJ (solid line with $\circ$), sHC (short dashed line with $\triangle$), BJ (dash-dot line with
$\circ$), HC (long dashed line with $\triangle$) and penalized scan (dotted line with $+$) 
in the sparse case $\alpha=0.2$, $\beta=0.48$.
\label{fig:Comparison-of-moderate-sparse}}
\end{figure}

We set $\alpha=0.2$
and $\beta=0.48$, so $\beta/(1-\alpha)=0.6$.
Figure \ref{fig:Comparison-of-moderate-sparse}
shows that 
sHC and sBJ have much higher power than HC and BJ, as predicted
by our theory. The penalized scan still does very well.

\subsection{Simulation result for dense case}

\begin{figure}
\includegraphics[height=0.43\textheight]{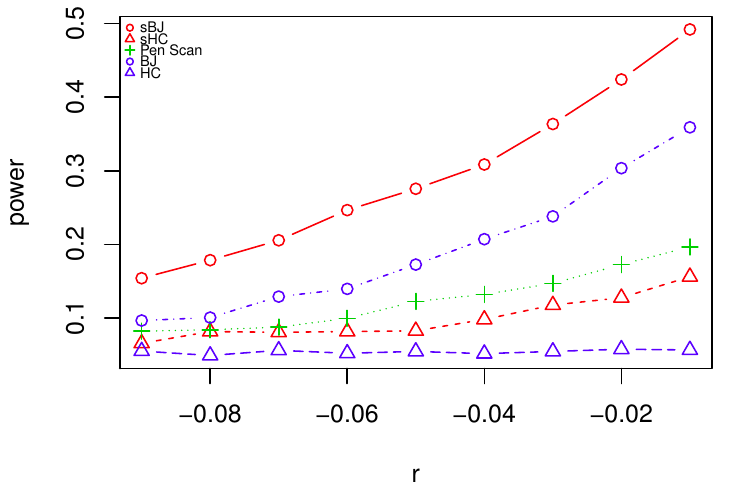}
\caption{Power for sBJ (solid line with $\circ$), sHC (short dashed line with $\triangle$), BJ (dash-dot line with
$\circ$), HC (long dashed line with $\triangle$) and penalized scan (dotted line with $+$)
in the dense case $\alpha=0.3$, $\beta=0.25$.
\label{fig:Comparison-of-Powers-dense}}
\end{figure}

We set $\alpha=0.3$
and $\beta=0.25$, so $\beta/(1-\alpha)=0.357$.
Since we are now in the dense regime, the scale for $r$ is with respect to
the calibration (\ref{eq: dense-calibrate}). 
While all five methods are asymptotically optimal in this situation,
Figure~\ref{fig:Comparison-of-Powers-dense} 
shows that there is quite some spread in the performance in this finite sample
setting. This reflects the observation in \cite{walther2013average}
that for these types of problems the asymptotics set in only slowly
and that performance should be assessed by simulations.
sBJ is the clear
winner in this case. HC and sHC have the worst performance, which
is the flip side of the effect described in 
Section~\ref{subsec:Simulation-result-for-very-sparse} as the relevant
information is now contained away from the smallest p-values.
Moreover, we can see that the structured versions of HC and BJ are
more powerful than their original counterparts, which indicates
sHC and sBJ can take some advantage of the structure in the signal
even in the dense case.

\section{The multivariate case\label{multivariate}}

All of the previous results can be readily extended to a multivariate
setting. We will use the superscript
$(d)$ to denote the dimension.
In order to keep the notation simple we will focus on the bivariate
case which already contains all the relevant ideas. 
The model (\ref{eq: Def-Block-mixture-problem}) then becomes

\begin{equation}
X_{ij}=\mu\boldsymbol{1}_{\cup_{g=1}^{m}I_{g}}(i,j)+Z_{ij},\qquad i,j=1,\ldots,n,
\label{eq: two-dimension-problem}
\end{equation}

\nin where the $Z_{ij}$ are i.i.d. standard normal and $\boldsymbol{1}_{\cup_{g=1}^{m}I_{g}}(i,j) =1$
iff the grid point $(i,j)$ is contained in an axis-parallel rectangle $I_g$ for
some $g \in 1,\ldots,m$. 
Analogously to the univariate case we assume that the rectangles $I_g$ are mutually disjoint
and randomly located on the Cartesian grid $\{1,\ldots,n\}^2$.
The number of axis-parallel rectangles (blocks) is now parametrized by
$m=n^{2(1-\alpha-\beta)}$ and each unknown 
rectangle $I_{g}$ contains $|I_{g}|=n^{2\alpha}$ grid points,
where $0\le\alpha<1$ and $0<\alpha+\beta\le1$. 

The task is to test $H_{0}:\mu=0$ vs.
$H_{1}:\mu=\mu(n)>0$. It was seen in the univariate case that the construction of
an appropriate approximating set is critical for optimally aggregrating the
information within and across blocks. This univariate approximating set can be easily
extended to the multivariate situation be taking cross-products:
Recall that in the univariate case the approximation set $\Iapp (\ell)$
depends on a precision parameter $\epsilon_{\ell}$. We now make this dependence
explicit by writing $\Iapp (\ell,\epsilon_{\ell})$ for this univariate collection. 
Now we construct a multivariate approximation 
set for axis-parallel rectangles in $\{1,\ldots,n\}^d$ via
the cross-product of univariate approximation sets $\bigtimes_{i=1}^d \Iapp (\ell_i,\eps_{\ell})$,
where the precision parameter $\eps_{\ell}$ depends on the volume of the rectangle but the
$\ell_i$ may vary to allow various aspect ratios:

For each level $\ell =0,\ldots, \lmax:= \lceil \log_2 (\frac{n}{8})^d \rceil$ we
set
$$
\Iapp ^{(d)} (\ell)\ :=\ \Bigl\{ R:\ 2^{\ell -1}< |R| \leq 2^{\ell} \mbox{ and }
R \in \bigtimes_{i=1}^d  \Iapp(\ell_i,\eps_{\ell}) \mbox{ for some }
0\leq \ell_i \leq \lceil \log_2 \frac{n}{8} \rceil \Bigr\}
$$
with $\eps_{\ell} =\frac{1}{6 \sqrt{\log_2 \frac{n^d}{2^{\ell -1}}}}$.
While this construction is somewhat different from that given
in \cite{walther2010optimal} for the density case, it enjoys similar
properties, see Proposition~\ref{properties} in Section~\ref{sec:Proofs}.
In particular, the cardinality of $\bigcup_{\ell} \Iapp ^{(d)} (\ell)$ is
$O\Bigl(n^d (\log n)^d\Bigr)$, so relevant computation can be done in time that is
almost linear in the number of observations $n^d$.

Now we can construct our test statistics exactly as in the univariate case: 
The \textit{structured higher criticism} $sHC_{n}^{(d)}$
and \textit{structured Berk-Jones} statistic $sBJ_{n}^{(d)}$ are defined as follows:
\[
sHC_{n}^{(d)}=\max_{\ell=0}^{\lmax}\sqrt{\frac{n^d}{2^{\ell}n_{\ell}}}HC_{n_{\ell}}(\ell),
\]

\[
sBJ_{n}^{(d)}=\max_{\ell=0}^{\lmax}\frac{n^d}{2^{\ell}n_{\ell}}BJ{}_{n_{\ell}}(\ell),
\]

\nin where $HC_{n_{\ell}}(\ell)$ and $BJ{}_{n_{\ell}}(\ell)$ denote the one-sided higher
criticism (\ref{HC}) and Berk-Jones statistic (\ref{BJ})
evaluated on the $n_{\ell}:=\# \Iapp ^{(d)} (\ell)$ p-values
pertaining to $\Iapp ^{(d)} (\ell)$, i.e. the p-values
$\{\bar{\Phi} (\sum_{(i,j)\in I}X_{ij}/\sqrt{|I|}), I\in \Iapp ^{(d)} (\ell)\}$.

Note that the definition of these structured statistics differs from the univariate
case only in the rescaling factor $\frac{n^{d}}{2^{\ell}}$ in place of $\frac{n}{2^{\ell}}$.
This is due to the fact that we now have an array of $n^d$ observations rather than $n$. Thus
there are now about 
$\frac{n^{d}}{2^{\ell}}$ disjoint intervals in $\Iapp ^{(d)} (\ell)$,
hence about $\frac{n^{d}}{2^{\ell}}$ of the p-values are independent.

We now focus on the bivariate case and establish the null distribution of these
statistics:

\begin{theorem}
\label{thm:multivariate-null}Under the null hypothesis $\mu=0$,
\begin{align*}
\frac{sBJ_{n}^{(2)}}{\log\log n}\  & \stackrel{p}{\leq}\ \frac{80}{9}\ (n\rightarrow\infty)\\
sHC_{n}^{(2)}\  & =\ O_{p}(\log^{4}n)
\end{align*}
\end{theorem}

The lower bound for detection in the sparse case $\frac{\beta}{1-\alpha}>\frac{1}{2}$
is the same as in the univariate setting after accounting for the sample size $n^2$ in place of $n$,
and $sHC_{n}^{(2)}$ and $sBJ_{n}^{(2)}$ are asymptotically optimal for detection:

\begin{theorem}
\label{optimality2}
The conclusions of Theorems~\ref{thm: Block-mixture-detection-boundary}(i) and
\ref{thm: Spatial HC is optimal}(i)
continue to hold for the model~(\ref{eq: two-dimension-problem}) with $\frac{\beta}{1-\alpha}>\frac{1}{2}$
and the calibration
$$
\mu=\mu(n)=\sqrt{2r\log n^2}/\sqrt{n^{2\alp}}.
$$
That is, if $r<\rho^{*}(\alpha,\beta)$, where $\rho^{*}$ is given in (\ref{eq: detection_boundary_sparse}),
then $H_{0}$ and $H_{1}$ merge asymptotically, i.e. the sum of
Type I and Type II errors tends to 1 for any test. If $r>\rho^{*}(\alpha,\beta)$,
then $sHC_{n}^{(2)}$ and $sBJ_{n}^{(2)}$ have asymptotic power 1.
\end{theorem}

\section{\label{graphs} Clusters in a network on the square lattice}

This section concerns the problem of detecting whether in a given network,
e.g. in a network of sensors, there are clusters of nodes that exhibit an
``unusual behavior''. This setting is important for a number of applications,
e.g. in surveillance, environmental monitoring and disease outbreak detection,
see \cite{ACCD} who treat the case of detecting a single (or a small number
of) clusters in a network.

Here we show how the evidence of such unusual behavior can be aggregated
over many such clusters. We follow \cite{ACCD} and model the network with the
$d$-dimensional square lattice. For simplicity we will derive our results
for the case $d=2$, which already contains all the essential ideas. We are
interested in the case where the signal is present on graph neighborhoods of
vertices, which we model as open balls $B_r(x)$ with center $x \in \{1,\ldots,n\}^2$
and radius $r$. The results in this section hold for balls with respect to the
$\ell^1$-norm, which corresponds to the shortest-path distance in a graph,
as well as the Euclidean norm. We derive our results for the latter as this is
the technically more demanding case, see Lemma~\ref{CL}.
Our model is therefore
\begin{equation}
X_{ij}=\mu\boldsymbol{1}_{\cup_{g=1}^{m}N_{g}}(i,j)+Z_{ij},\qquad i,j=1,\ldots,n,
\label{C*}
\end{equation}
where the $Z_{ij}$ are i.i.d. standard normal and each graph neighborhood $N_g$
is a ball with respect to the $\ell^2$-norm (or the $\ell^1$-norm) that contains
$|N_g|=n^{2 \alp}$ grid points, where $0\leq \alp <1$. As before, we assume that the
$N_g$ are mutually disjoint and randomly located on the Cartesian grid $\{1,\ldots,n\}^2$
and the number of balls is parametrized by $m=n^{2(1-\alp -\beta)}$.

The task is to test $H_0:\ \mu=0$ vs. $H_1:\ \mu=\mu(n)>0$. In order to apply
the general recipe of this paper for optimally aggregating the information within
and across neighborhoods, we need to construct an appropriate approximating set
for the neighborhoods. The idea for this construction can be adapted from the
previous settings, which shows the generality of this approach:

We approximate balls with volume in $(\pi 2^{\ell -1},\pi 2^{\ell}]$, where $\ell=0,\ldots,
\lmax=\lceil \log_2 \frac{n^2}{8} \rceil$, with the collection
$$
\Capp (\ell) := \Bigl\{ B_{r_i}(j,k): r_i^2=2^{\ell -1 +i \eps_{\ell}},\, 
i=0,\ldots, \Bigl\lfloor \frac{1}{\eps_{\ell}} \Bigr\rfloor,\ 
j,k \in \{m\, d_{\ell}, m \in \N\} \cap [r_i, n-r_i+1]\Br \}
$$
where $\eps_{\ell}:=\frac{1}{\sqrt{\log_2 \frac{n^2}{2^{\ell-1}}}}$, $d_{\ell}
=\lceil \eps_{\ell} 2^{\frac{\ell -1}{2}} \rceil$.
That is, we approximate the centers with a grid whose spacing is a fraction
$\eps_{\ell}$ of the square root of the approximate volume of the ball, $2^{\ell -1}$,
and we approximate the square radius with a geometric progression. Proposition~\ref{PropC}
shows that the balls in $\bigcup_{\ell} \Capp (\ell)$ can approximate every ball
with small relative error, while the cardinality of $\bigcup_{\ell} \Capp (\ell)$
is almost linear in the sample size $n^2$:
\begin{proposition} \label{PropC}

\noindent \begin{itemize}
\item[(i)] $\# \bigcup_{\ell} \Capp (\ell)\ =\ O\Bigl( n^2 (\log n)^{\frac{3}{2}} \Bigr)$ \\
\item[(ii)] For every ball $B_R(s,t)$ with $R^2 \in [1,\frac{n^2}{8}]$ and $s,t \in [R,n-R+1]$
there exists $B_r(j,k) \in \bigcup_{\ell} \Capp (\ell)$ such that
$$
| B_R(s,t) \triangle B_r(j,k) | \ \leq \ 3 \frac{|B_R(s,t)|}{\sqrt{ \log_2 \frac{n^2}{|B_R(s,t)|}}}.
$$
\end{itemize}
\end{proposition}

Furthermore, it will be shown in the proof of Theorem~\ref{CD} that the balls
in $\Capp (\ell)$ can be grouped into a small number of
at most $8 (\log n)^{1/2}$ of groups such that each group contains $\sim \frac{n^2}{2^{\ell +2}}$ 
mutually disjoint balls. This allows to define the structured higher criticism and Berk-Jones
statistics as in Section~\ref{multivariate},
where as before $HC_{n_{\ell}}(\ell)$ and $BJ{}_{n_{\ell}}(\ell)$ denote the one-sided higher
criticism (\ref{HC}) and Berk-Jones statistic (\ref{BJ})
evaluated on the $n_{\ell}:=\# \Capp (\ell)$ p-values
pertaining to $\Capp (\ell)$, i.e. the p-values 
$\{\bar{\Phi} (\sum_{(i,j)\in I}X_{ij}/\sqrt{|I|}), I\in \Capp (\ell)\}$.
As a consequence we obtain results for the null distributions of these statistics and
optimality properties that are analogous to those for univariate and multivariate rectangles:

\begin{theorem}
\label{CD}
Under the null hypothesis $\mu=0$ there exists $C>0$ such that
\begin{align*}
\frac{sBJ_{n}^{(2)}}{\log\log n}\  & \stackrel{p}{\leq}\ C\ (n\rightarrow\infty)\\
sHC_{n}^{(2)}\  & =\ O_{p}(\log^{\frac{5}{2}}n).
\end{align*}
Moreover, the conclusions of Theorem~\ref{optimality2} continue to hold for model (\ref{C*}).
\end{theorem}

\section{Composite alternatives\label{composite}}

We developed the above theory for the one-sided alternative
$H_{1}:\mu=\mu(n)>0$ in the models (\ref{eq: Def-Block-mixture-problem},
\ref{eq: two-dimension-problem},\ref{C*}).
A more general alternative obtains by allowing the sign of $\mu$ to change from
block to block. Then model (\ref{eq: Def-Block-mixture-problem}) becomes
\be
X_{i}=\sum_{g=1}^m (-1)^{s_g}\mu\boldsymbol{1}_{I_{g}}(i)+Z_{i},\qquad i=1,\ldots,n,
\ee
where $s_g \in \{1,2\}$. Such alternatives can be tested by applying the structured
test statistic to the p-values from two-sided z-tests, as suggested by
\cite{delaigle2009higher} in the context of the higher criticism statistic. That is,
one employs the two-sided p-values $2\bar{\Phi} (|\X(I)|)$ in place of $\bar{\Phi} (\X(I))$.
This test procedure satisfies the same
optimality results that we derived for the one-sided case since the proofs of these results
depend on establishing certain polynomial growth rates up to logarithmic factors, while
the use of two-sided p-values affects these rates only with a factor 2.
Likewise, the optimality results for the multivariate and network
settings (\ref{eq: two-dimension-problem},\ref{C*}) continue to hold for two-sided p-values.

An alternative approach is to base the analysis not on scaled averages $\X(I)$ but on 
mean squares. The technical analysis of the resulting structured test statistic
will be somewhat different since it involves the tails of noncentral chi-squared distributions
rather than Gaussians. In the context of the sparse signal model,
\cite{donoho2004higher} find that the higher criticism statistic based on data
from a chi-squared distribution achieves the same optimal detection region as
in the Gaussian case.

\section{Discussion\label{sec:Conclusion-and-Discuss}}

In this paper, we established the lower bound for detection in the
multiple blocks model. An asymptotically optimal method is also
proposed which is adaptive to the unknown number of blocks and to the unknown
block length. It was shown how this methodology can be readily extended to the
multivariate situation and to detecting clusters in a network.

Another interesting problem for future research is the identification version of this
problem, in which we not only want to detect whether a signal is present,
but we also want to approximately find the location of all blocks of
signals. In \cite{kou2017identifying} it is shown that when
there is only one block of signals (corresponding to $\alpha +\beta=1)$, then
the identification and the detection problem are of the same difficulty.
However, in the more general case where $\alpha +\beta<1$ some calculations
show that the identification problem is necessarily more difficult
than the detection problem, in the sense that the lower bound for
the former is larger. To the best of our knowledge, an adaptively optimal method
is not yet known for the corresponding multiple blocks
identification problem. We leave this as an open problem for future research.

\section{Proofs\label{sec:Proofs}}

\subsection{Some basic results}

It is helpful to analyze the statistical behavior of the test statistics via tail
probabilities. To this end, note that since $\bar{\Phi}$ is strictly decreasing we
have the following representation for distinct real numbers $X_1,\ldots,X_n$ and
$p_i :=\bar{\Phi} (X_i)$:
\be  \label{represent}
\frac{k}{n}-p_{(k)}\ =\ \frac{1}{n} \sum_{i=1}^n \Bl( {\bf 1}(X_i \geq
\bar{\Phi}^{-1} (p_{(k)})) -p_{(k)}\Br),\ \ \ \ k=1,\ldots,n.
\ee 

The following Lemma summarizes important properties of the univariate
approximating set $\Iapp (\ell)$ while Proposition~\ref{properties} gives
some relevant results for the multivariate version
$ \Iapp^{(d)}(\ell)$ in dimensions $d \geq 1$. Further results for the bivariate
case can be found in Lemma~\ref{lemmamult}.

\begin{lemma}
\label{lem: Short-range-dependent}
The intervals in $\Iapp(\ell)$
can be grouped into at most $\min (2^{2\ell},4 \eps_{\ell}^{-2}) \leq
144 \log_2 n$ groups such
that each group consists of either $\lfloor\frac{n}{L_{\ell}}\rfloor$
or \textup{$\lfloor\frac{n}{L_{\ell}}\rfloor-1\ge\lfloor\frac{n}{2^{\ell}}\rfloor-1$}
disjoint intervals,
where $L_{\ell}$ is the largest multiple of $d_{\ell}$ that is not
larger than $2^{\ell}$. Further,
$\#\Iapp (\ell)\leq n2^{-\ell} \min (2^{2\ell},4 \eps_{\ell}^{-2})
\leq 144 n 2^{-\ell} \log_2 n$.
\end{lemma}

\begin{proof}[Proof of Lemma \ref{lem: Short-range-dependent}]
Let $S_{\ell}$ be the collection of all intervals in $\Iapp(\ell)$ whose left endpoint is smaller
than $L_{\ell}$, where $L_{\ell}$ is the largest multiple of $d_{\ell}$ that is not larger than
$2^{\ell}$. For a given $I \in S_{\ell}$ consider the collection of shifts of $I$ by multiples of
$L_{\ell}$: $\mbox{shift}_{\ell}(I):= \Bl\{ J \subset (0,n]:\, J=kL_{\ell} +I,\ k=0,1,2,\ldots\Br\}$.
Since $I \in \Iapp(\ell)$ implies $|I|\leq L_{\ell}$, the intervals in $\mbox{shift}_{\ell}(I)$ are
disjoint and there are either $\lfloor \frac{n}{L_{\ell}} \rfloor$ or
$\lfloor \frac{n}{L_{\ell}} \rfloor -1 \geq \lfloor\frac{n}{2^{\ell}}\rfloor -1$ intervals in
$\mbox{shift}_{\ell}(I)$. One readily observes that each interval $I \in \Iapp(\ell)$ can be generated
by such a shift: $$ \Iapp(\ell) \ =\ \bigcup_{I \in S_{\ell}} \mbox{shift}_{\ell}(I). $$ Finally, there
are exactly $\frac{L_{\ell}}{d_{\ell}}$ different starting points for intervals $I \in S_{\ell}$.
Since each such interval $I$ satisfies $2^{\ell -1} < |I| \leq 2^{\ell}$ we obtain for $\ell \geq 1$:
$$ \# S_{\ell}\ \leq \ \frac{L_{\ell}}{d_{\ell}} \Bl( \frac{2^{\ell}-2^{\ell -1}}{
d_{\ell}} +1 \Br)\ \leq \
\Bl(\frac{2^{\ell}}{d_{\ell}}\Br)^2 \ \leq \
\min (2^{2\ell},4 \eps_{\ell}^{-2}) \ \leq \ 144 \log_2 n
$$
and the same bound holds for $\ell =0$. As for the upper bound on $\#\Iapp(\ell)$,
an analogous counting argument shows that there are not more than $\frac{n}{d_{\ell}}$ starting
points, each having not more than $\lceil\frac{2^{\ell}-2^{\ell-1}}{d_{\ell}}\rceil
\leq \frac{2^{\ell}}{d_{\ell}}$ endpoints if $\ell \geq 1$. Hence
$\#\Iapp(\ell) \leq \frac{n2^{\ell}}{d_{\ell}^2}$ and the claimed bound follows
from the above inequality; the same bound clearly also holds for $\ell =0$.
\end{proof}
\smallskip

\begin{proposition}  \label{properties}

\noindent \begin{itemize}
\item[(i)] $\# \bigcup_{\ell} \Iapp ^{(d)} (\ell)\ =\ O\Bigl(n^d (\log n)^d\Bigr)$
\item[(ii)] For every axis-parallel rectangle $R \subset \{1,\ldots,n\}^d$ with sides 
not longer than $\frac{n}{8}$ there exists
$\tilde{R} \in \bigcup_{\ell} \Iapp ^{(d)} (\ell)$ such that $|R \triangle \tilde{R}|
\leq C_d \frac{|R|}{\sqrt{\log_2 \frac{n^d}{|R|}}}$ for some universal constant $C_d$.
\item[(iii)] The definition of $\Iapp ^{(d)} (\ell)$ implies the constraint $\ell \leq
\sum_{i=1}^d \ell_i \leq \ell +d-1$ for the marginal levels $\ell_i$.
\end{itemize}
\end{proposition}

Employing the latter constraint is helpful for efficiently enumerating the rectangles
in $\Iapp ^{(d)} (\ell)$, e.g. in simulations. As an aside, if one modifies the definition
of $\Iapp ^{(d)} (\ell)$ to let the $\ell_i$ be as large as $\lceil \log_2 n \rceil$
and $\ell$ as large as $\log_2 (n^d/8)$, then $\Iapp ^{(d)} (\ell)$ will also contain
approximating rectangles for all marginal distributions.

\begin{proof}[Proof of Proposition~\ref{properties}] 
Lemma \ref{lem: Short-range-dependent} gives $\# \Iapp(\ell_i,\eps_{\ell})
\leq n 2^{-\ell_i} \left(4\eps_{\ell}^{-2}\right)
\leq$ \linebreak[4] $144 n 2^{-\ell_i} \log_2 n^d$. Hence $\# \bigcup_{\ell_i =0}^{\lceil \log_2
\frac{n}{8} \rceil} \Iapp (\ell_i, \eps_{\ell}) \leq 288 d n \log_2 n$ and so
$$
\# \bigcup_{\ell} \Iapp^{(d)} (\ell)\ \leq \ \# \bigtimes_{i=1}^d
\bigcup_{\ell_i =0}^{\lceil \log_2
\frac{n}{8} \rceil} \Iapp (\ell_i, \eps_{\ell}) \ \leq \ (288 dn \log_2 n)^d
$$
proving (i).

As for (ii), let $R=I_1 \times \ldots \times I_d$ be an axis-parallel rectangle,
so each $I_i$ is an interval of the form $(j_i,k_i]\subset (0,n]$ with length at
most $n/8$. Hence there exists $\ell_i \in \{0,\ldots,\lceil \log_2 \frac{n}{8} \rceil \}$
such that $2^{\ell_i-1}< |I_i| \leq 2^{\ell_i}$, and there exists 
$\ell \in \{0,\ldots,\lceil \log_2 (\frac{n}{8})^d \rceil \}$ such that
$2^{\ell-1}< |R| \leq 2^{\ell}$. So by the definition of $\Iapp (\ell_i,\eps_{\ell})$,
there exists $\tilde{I}_i \in \Iapp (\ell_i,\eps_{\ell})$ with $|I_i \triangle \tilde{I}_i|
\leq 2d_{\ell_i} \leq 2 \eps_{\ell} |I_i|$ for $i=1,\ldots,d$. Thus, by decomposing
$R \triangle \tilde{I}$ and collecting terms, we get $|R \triangle \tilde{I}|
\leq C_d \eps_{\ell} |R|$ for a constant $C_d$. Finally, $\eps_{\ell} =
\frac{1}{6\sqrt{\log_2 \frac{n^d}{2^{\ell-1}}}} \leq \frac{1}{6\sqrt{\log_2 \frac{n^d}{|R|}}}$.
(If we arrange $\tilde{I}_i \subset I_i$ for all $i$ by modifying the definition of
$\Iapp$ somewhat, then clearly $|R \triangle \tilde{I}|\leq 2d \eps_{\ell} |R|$.) 

Concerning (iii), $R \in \Iapp ^{(d)}(\ell)$ implies $R=I_1 \times \ldots \times I_d
\in \bigtimes_{i=1}^d  \Iapp(\ell_i,\eps_{\ell})$. So $2^{\ell_i-1}<|I_i| \leq 2^{\ell_i}$
and $2^{\ell -1} < \prod_{i=1}^d |I_i| \leq 2^{\ell}$, hence $\ell -1 < \sum_i \ell_i$
and $\sum_i (\ell_i -1) < \ell$. 
\end{proof}

\subsection{Proofs for Section~\ref{sec:Def-Block-heterogeneous-mixtures}}

\begin{proof}[Proof of Theorem \ref{thm: Block-mixture-detection-boundary}](i)
We may assume without loss of generality that $\frac{n}{n^{\alp}}$ is an integer.
Denote by (A) the submodel where the signals can only start and end on a grid given
by $\{i (n^{\alp})+1,\ldots(i+1)n^{\alp}\}$ for $i=0,\ldots,\frac{n}{n^{\alp}}-1$.
It is enough to show that Theorem
\ref{thm: Block-mixture-detection-boundary} holds for this submodel (A) since detection
in the submodel is not more difficult than in the original model and hence the detection
boundary for the submodel cannot be larger than for the original model.

Let $S_{i}:=\sum_{j=1}^{n^{\alp}}X_{(i-1)n^{\alp}+j}/\sqrt{n^{\alp}}=:s_{i}+Z_{i}^{'}$
for $i=1,\ldots,n^{'}$, where $n^{'}=\frac{n}{n^{\alp}}=n^{1-\alpha}$.
Then $Z_{i}^{'}\stackrel{iid}{\sim}N(0,1)$, and $s_{i}=0$ for all
but $n^{1-\alpha-\beta}$ locations, while at these locations $s_{i}=\sqrt{2r\log n}=\sqrt{2r^{'}\log n^{'}}$,
where $r^{'}=\frac{r}{1-\alpha}$. The locations of these elevated means
are a random sample without replacement from
$\{1,\ldots,n^{'}\}$. We denote this
model by (B). (B) is in fact a sparse signal model (\ref{HM})
with $n^{'}=n^{1-\alpha}$,
and sparsity level $\beta^{'}=\frac{\beta}{1-\alpha}>\frac{1}{2}$. It was 
proved in \cite{Ingster97,ingster1998minimax,IngsterSuslina}, 
see also section~1.1 in \cite{donoho2004higher}, that the lower bound 
for model (B) is given by (\ref{HMboundary}) with $\beta^{'}$ in place
of $\beta$. (In some of these proofs the number of elevated means follows a binomial distribution,
but the result is regularly referenced for the setting (\ref{HM}) where that number is fixed, 
see e.g. \cite{hall2010innovated,zhong2013tests}.
Theorem~4 in \cite{Ingster97} shows that the result does indeed carry over to the setting (\ref{HM}),
and the proof in \cite{IngsterSuslina} is explicitly for this setting.)

Writing $\mu^{'}:=\sqrt{2r \log n} =\mu \sqrt{n^{\alp}}$ and observing 
$$
-\sum_{j=1}^{n^{\alp}}\frac{\left( X_{(i-1)n^{\alp}+j} -\mu\right)^2}{2}
+ \sum_{j=1}^{n^{\alp}} \frac{X_{(i-1)n^{\alp}+j}^2}{2} \ =\ S_i \sqrt{n^{\alp}} \mu -
\frac{\left( \sqrt{n^{\alp}} \mu \right)^2}{2} \ =\ -\frac{ (S_i -\mu^{'} )^2}{2} +\frac{S_i^2}{2}
$$
shows that the likelihood ratio test has the same
test result on model (A) and (B). Therefore, written in our original notation
$\alpha,\beta,r$, the lower bound for model (A) gives 
(\ref{eq: detection_boundary_sparse}). 

The proof of (ii) is analogous to (i),
but now $\beta^{'}=\frac{\beta}{1-\alpha}<\frac{1}{2}$ and
the elevated means are $s_{i}=n^{-r}=(n^{'})^{-r^{'}}$.
The lower bound (\ref{HMdenseBoundary}) established in \cite[Theorem~3]{tony2011optimal}
translates into (\ref{eq: detection_boundary_dense}).
\end{proof}

\subsection{Proofs for Section~\ref{sec:Spatial-HC-BJ-null}}

\begin{proof}[Proof of Theorem \ref{thm:HC-BJ-Under-the-null}]
To prove the theorem, we fix $\ell$ and apply Lemma~\ref{lem: Short-range-dependent}.
Write $G_i$, $i=1,\ldots,i_{\mathrm{max}}$,
  for the set of p-values pertaining to the intervals in the $i$th group given by  
Lemma~\ref{lem: Short-range-dependent}. As those intervals are disjoint, these p-values are i.i.d. 
$U[0,1]$ under $H_0$. Further $\sum_{i=1}^{i_{\mathrm{max}}} \# G_i = n_{\ell}$. Denote by $F^{(i)}$ the 
empirical cdf of the p-values in $G_i$. Then we can write the empirical cdf of all $n_{\ell}$ 
p-values as $F_{n_{\ell}} = \sum_{i=1}^{i_{\mathrm{max}}} \frac{\# G_i}{n_{\ell}} F^{(i)}$. 
Recall that $BJ_{n_{\ell}}(\ell)$ is defined by (\ref{BJ}) evaluated at these
$n_{\ell}$ p-values. Hence 
$$ 
BJ_{n_{\ell}}(\ell) \leq \sup_{t \in [p_{(1)},p_{(n_{\ell})}]} 
n_{\ell}  \Bigl( F_{n_{\ell}}(t) \log \frac{F_{n_{\ell}}(t)}{t}   +(1-F_{n_{\ell}}(t)) \log 
\frac{1-F_{n_{\ell}}(t)}{1-t} \Br). 
$$ 
Since the function $(s,t) \ra s\log \frac{s}{t} +(1-s) \log 
\frac{ 1-s}{1-t}$ is convex on $(0,1)^2$, Jensen's inequality gives 

\begin{eqnarray*} 
&  & \frac{n}{2^{\ell}n_{\ell}} BJ_{n_{\ell}}(\ell)   \\
& \leq & \frac{n}{2^{\ell}n_{\ell}} \sup_{t \in [p_{(1)},p_{(n_{\ell})}]}     
n_{\ell} \sum_{i=1}^{i_{\mathrm{max}}} \frac{\# G_i}{n_{\ell}}     
\Bigl( F^{(i)}(t) \log \frac{F^{(i)}(t)}{t}     
+(1-F^{(i)}(t)) \log \frac{1-F^{(i)}(t)}{1-t} \Br)\\ 
& \leq & \max_{i=1,\ldots, i_{\mathrm{max}}} \sup_{t \in [p_{(1)},p_{(n_{\ell})}]}   
\frac{n}{2^{\ell}} \Bigl( F^{(i)}(t) \log \frac{F^{(i)}(t)}{t}   
+(1-F^{(i)}(t)) \log \frac{1-F^{(i)}(t)}{1-t} \Br). 
\end{eqnarray*}

The last inequality is conservative as we bound the weighted average of $i_{\mathrm{max}}$ 
Berk-Jones statistics by the worst case; obtaining a better bound is not straightforward as the 
Berk-Jones statistics are dependent.
Setting $A:=p_{(1)}$ and $B:=p_{(n_{\ell})}$ in the proof of  of the third inequality of 
Proposition~\ref{prop:HC-BJ-each-lvel-null} shows that for every $\eta>0$, $K>1$, and for
every group $i$: 

\begin{eqnarray*} 
& & \P_{H_0} \Bl( \sup_{t \in [p_{(1)},p_{(n_{\ell})}]}  \# G_i \Bl(   
F^{(i)}(t) \log \frac{F^{(i)}(t)}{t}   +(1-F^{(i)}(t)) \log \frac{1-F^{(i)}(t)}{1-t} \Br) > \eta \Br)\\ 
& \leq & 22K (\log \# G_i) (\eta +1) \exp (-\eta) +   
\P_{H_0} \Bl( p_{(1)} < \frac{1}{n^K} \mbox{ or } p_{(n_{\ell})} >   1 -\frac{1}{n^K}\Br) \\
& \leq & 22K (\log n) (\eta +1) \exp (-\eta) + 288 \frac{(\log_2 n)n}{n^K} \end{eqnarray*}
since Lemma~\ref{lem: Short-range-dependent} gives $\# G_i \leq n$ and 
$n_{\ell} \leq 144 n\log_2 n $. Further Lemma~\ref{lem: Short-range-dependent} gives 
$\lfloor \frac{n}{2^{\ell}}\rfloor \leq \# G_i +1$ for all $i$. 
For simplicity of exposition we will use $\frac{n}{2^{\ell}} \leq \# G_i$ (the remainder of the 
proof can be readily adapted to the weaker condition with  standard arguments). 
Applying the union bound first over $i \leq i_{\mathrm{max}}$ (and noting $i_{\mathrm{max}} \leq 144\log_2 n$ 
by Lemma~\ref{lem: Short-range-dependent})
 and then over $\ell \leq \lmax$ gives for $\eta =c \log \log n$: 
\begin{eqnarray*}
& & \P_{H_0} \Bl( sBJ_n^{} > c \log \log n \Br)\ \\
& \leq & (\lmax+1) 144 \log_2 n \Bl[ 22K (\log n)^{1-c} 
(c \log \log n +1)+288 \frac{(\log_2 n)n}{n^K} \Br]. 
\end{eqnarray*}
Since $\lmax \leq \log_2 n$ this bound will converge to 0 for  $c>3$ and $K>1$, 
proving the claim for $sBJ_n^{}$.      
Concerning $sHC_n^{}$, as in (\ref{HCst}) we get 
\begin{align*} 
\sqrt{\frac{n}{2^{\ell}n_{\ell}}} HC_{n_{\ell}}(\ell) & \leq \sqrt{\frac{n}{2^{\ell}n_{\ell}}} 
\sup_{t \in [p_{(1)},p_{(n_{\ell})}]}   \sqrt{n_{\ell}} \,\frac{F_{n_\ell}(t) -t}{\sqrt{t(1-t)}} \\ 
& \leq \sqrt{\frac{n}{2^{\ell}}} \sup_{t \in [p_{(1)},p_{(n_{\ell})}]}   
\sum_{i=1}^{i_{\mathrm{max}}} \frac{\# G_i}{n_{\ell}} \,\frac{F^{(i)}(t)-t}{   
\sqrt{t(1-t)}} \\ 
& \leq \max_{i=1,\ldots, i_{\mathrm{max}}} \sup_{t \in (0,1)} 
\sqrt{\frac{n}{2^{\ell}}}   \,\frac{F^{(i)}(t)-t}{\sqrt{t(1-t)}} 
\end{align*}
Using $\frac{n}{2^{\ell}} \leq \# G_i$ for all $i$, the last inequality of 
Proposition~\ref{prop:HC-BJ-each-lvel-null}, and applying the union bound over  
$i \leq i_{\mathrm{max}}$ (noting $i_{\mathrm{max}} \leq 
144 \log_2 n$) and $\ell \leq \lmax$ then gives for $\eta =B \log^2 n$ with $B \geq 1$: 
$$ \P_{H_0} \Bl( sHC_n^{} > B \log^2 n \Br)\ \leq \ \frac{C(\lmax+1) 144 \log_2 n}{B 
\log^2 n} 
$$ 
for some $C$ not depending on $B$. The claim follows as $\lmax \leq \log_2 n$. 
\end{proof}

\begin{proof}[Proof of Proposition \ref{prop:HC-BJ-each-lvel-null}]
It is a well known fact that $B(t)=(1+t)U(t/(1+t))$ is a standard Brownian
motion, for which \cite[p.34]{ito1965diffusion} establish the
following inequality:
\[
\P\Bigl(\sup_{t\in[a,b]}\frac{B(t)}{f(t)}>1\Bigr)\le
\int_{0}^{a/f(a)^{2}}\frac{e^{-1/(2t)}}{\sqrt{2\pi t^{3}}}\, dt
+\int_{a}^{b}\frac{f(t)}{\sqrt{2\pi t^{3}}}e^{-f(t)^{2}/(2t)}dt
\]
for $0<a<b\le1$ and $f(t)$ increasing on $(0,b]$. Setting
$f(t)=\eta\sqrt{t}$ we obtain
\begin{align*}
\P\Bigl(\sup_{t\in[a,b]}\frac{U(t)}{\sqrt{t(1-t)}}>\eta\Bigr) &=
\P\Bigl(\sup_{t\in[\frac{a}{1-a},\frac{b}{1-b}]}\frac{B(t)}{\sqrt{t}}>\eta
\Bigr) \\
 & =\int_{0}^{1/\eta^{2}}\frac{e^{-1/(2t)}}{\sqrt{2\pi t^{3}}}\,dt+
\int_{a/(1-a)}^{b/(1-b)}\frac{\eta}{\sqrt{2\pi}\,t}e^{-\eta^{2}/2}dt\\
 & = \int_{\eta}^{\infty} \frac{2}{\sqrt{2\pi}} e^{-t^2/2}dt+
\frac{\log\frac{b(1-a)}{a(1-b)}}{\sqrt{2\pi}}\,\eta \,e^{-\eta^{2}/2}\\
 & \le \frac{2}{\sqrt{2\pi}\,\eta} e^{-\eta^{2}/2}+
\frac{\log\frac{b(1-a)}{a(1-b)}}{\sqrt{2\pi}}\,\eta \, e^{-\eta^{2}/2} \\
 & = \frac{\frac{2}{\eta}+\eta \log\frac{b(1-a)}{a(1-b)}}{\sqrt{2\pi}}e^{-\eta^2/2}
\end{align*}
where we used Mill's ratio to bound the normal tail in the fourth line.

As for the second inequality,
Lemma~3.1 in \cite{duembgen2014confidence} gives for real $u$ and $c>0$: 
$$ \P \Bl( \sup_{t \in [l(u),l(u+c)]} n\left(F_{n}(t)\frac{\log F_{n}(t)}{t}+
(1-F_{n}(t)) \log\frac{1-F_{n}(t)}{1-t}\right) > \eta \Br)
 \ \leq \ 2 \exp \Bl( -e^{-c} \eta \Br), 
$$ 
where $l(u):=\frac{e^u}{1+e^u}$. Set $u=\log \frac{a}{1-a}$ and 
$c=\log \frac{b (1-a)}{a (1-b)} >0$, so $a=l(u)$, $b=l(u+c)$. Hence for any 
positive integer $K$: 
\begin{align*} 
\P \Bl( & \sup_{t \in [a,b]}  n\left(F_{n}(t)\frac{\log F_{n}(t)}{t}+
(1-F_{n}(t)) \log\frac{1-F_{n}(t)}{1-t}\right)> \eta \Br) \\
 & \leq \sum_{i=1}^K \P \Bl( \sup_{t \in [l(u+\frac{i-1}{K}c),   l(u+\frac{i}{K}c)]} 
n\left(F_{n}(t)\frac{\log F_{n}(t)}{t}+
(1-F_{n}(t)) \log\frac{1-F_{n}(t)}{1-t}\right)> \eta \Br) \\ 
& \leq 2K \exp \Bl( - e^{-\frac{c}{K}} \eta \Br) \\ 
& \leq 2\exp \Bl(- (1-\frac{c}{K}) \eta + \log K \Br) 
\end{align*} 
With a view towards minimizing this expression we set $K:=\lceil c \eta \rceil$. 
Then the above expression is not larger than 
$$ 
2\exp \Bl( - \eta +1 +\log \Bl\lceil c\eta \Br\rceil \Br) \ \leq \  
2e \Bl( \eta c+1 \Br) \exp \Bl(- \eta \Br). 
$$ 
As for the third inequality, elementary considerations show 
\begin{eqnarray} 
BJ_n\ & = & \ \sup_{t \in [U_{(1)},U_{(n/2)}]} n \Bigl( F_n(t) \log \frac{F_n(t)}{t}   +(1-F_n(t)) 
\log \frac{1-F_n(t)}{1-t} \Br)\ 1\Bl(t<F_n(t)\Br) \nn \\
& \leq & \sup_{t \in [A,B]} n \Bigl( F_n(t) \log \frac{F_n(t)}{t}   +(1-F_n(t)) \log 
\frac{1-F_n(t)}{1-t} \Br), \label{PropStar} 
\end{eqnarray} 
where $A=U_{(1)}, B=U_{(n)}$. For later reference it is convenient to prove the 
inequality for the latter statistic, i.e. the two-sided version of the Berk-Jones statistic that is 
based on all $n$ p-values rather than a fraction of them, and with general random limits 
$0 \leq A <B \leq 1$ for $t$.
For ease of notation let $K>1$ be such that $K \log_2 n$ is an integer. 
We will use the partition $[\frac{1}{n^K},\frac{1}{2}] = \bigcup_{i=1}^{K 
\log_2 n -1} [(\frac{1}{2})^{i+1},(\frac{1}{2})^i]$. Note that for each set in 
this partition we can apply the second inequality of the Proposition with the 
same exponential tail bound as the ratio of the right to the left endpoint is 2,
hence $\log \frac{b(1-a)}{a(1-b)} \leq \log 4$ as $b\leq \frac{1}{2}$.
 We can proceed analogously on $[\frac{1}{2}, 1-\frac{1}{n^K}]$ as the distribution 
of the statistic is symmetric about $\frac{1}{2}$. Applying the union bound to the 
resulting partition of $[\frac{1}{n^K}, 1-\frac{1}{n^K}]$ gives 
\begin{eqnarray*} 
& & \P_{H_0} \Bl( BJ_n > \eta \Br) \\
& \leq & 2(K \log_2 n -1) 2e (\eta \log 4 +1) \exp(-\eta)   \ + 
\ \P \Bl(A < \frac{1}{n^K} \ \mbox{ or } B > 1-\frac{1}{n^K} \Br) \\ 
& \leq & 22 K  (\log n) (\eta +1) \exp(-\eta) \ +\   \P \Bl(A < \frac{1}{n^K} \ \mbox{ or } 
B > 1-\frac{1}{n^K} \Br) . 
\end{eqnarray*} 
For $A=U_{(1)}, B=U_{(n)}$ the latter probability is not larger than $2n^{1-K}$, 
proving the claim for $BJ_n$.

Finally, elementary considerations show 
\be 
HC_n\ =\ \sup_{t \in [p_{(1)},p_{(n/2)}]} \sqrt{n} \frac{F_n(t)- t}{\sqrt{t(1-t)}}. 
\label{HCst} 
\ee
\cite[pp. 601--603]{shorackW} analyze $\sup_{t \in (0,1)} Z_n(t)$, 
where $Z_n(t)=\sqrt{n} \frac{F_n(t)-t}{\sqrt{t(1-t)}}$, by splitting $(0,1)$ into 
$[0,\frac{1}{n}]$, $[\frac{1}{n},d_n]$, $[d_n,\frac{1}{2}]$ (and their reflections 
about $\frac{1}{2}$), where $d_n=\frac{\log ^5 n}{n}$. The inequality they use for 
the first interval gives 
$$ 
\P \Bl( \sup_{t \in [0,\frac{1}{n}]} Z_n(t) > \eta \Br)\ \leq \ \Bl( 
\frac{\eta}{2} -1 \Br) ^{-1} \ \leq \ \frac{4}{\eta}\ \ \ \ \ \mbox{ for 
$\eta \geq 4$} 
$$ 
while the Shorack and Wellner inequality gives for the second interval 
\begin{align*} 
\P \Bl( \sup_{t \in [\frac{1}{n},d_n]} Z_n(t) > \eta \Br) & \leq  60(\log \log n) 
\exp\Bl(-\frac{3}{32} \eta \Br)\ \ \ \mbox{ if $\eta >  \frac{3}{2}$}\\  
& \leq \frac{C^{\prime}}{\eta}  \ \ \ \mbox{ for $\eta \geq \sqrt{\log \log n}$}. 
\end{align*}
On the interval $[d_n,\frac{1}{2}]$ one can use exponential inequalities for the 
Hungarian construction, see \cite[ch. 12.1]{shorackW}, as well 
as for $\frac{U(t)}{\sqrt{t(1-t)}}$, see above. The first shows that 
$\sup_{t \in [d_n,\frac{1}{2}]} |Z_n(t) -\frac{U(t)}{\sqrt{t(1-t)}}|$ satisfies the 
claimed tail bound whenever $\eta$ exceeds a certain constant, while the second gives the tail bound 
\begin{eqnarray*} 
& & \Bl(\frac{2}{\eta}+\eta \log 2 \frac{\frac{1}{2}}{d_n}\Br)
 \exp \Bl( -\frac{1}{2} \eta^2 \Br) \ \\
& \leq & \ \eta (\log n) \exp \Bl(-\frac{1}{D} \eta^2 \Br) \exp \Bl(-(\frac{1}{2} 
-\frac{1}{D}) \eta ^2 \Br)\\
& \leq & C^{\prime \prime} \eta (\log n) \exp \Bl(-\frac{1}{D} \eta^2 \Br)  
\eta^{-2} \ \ \ \mbox{ for some $C^{\prime \prime}=C^{\prime \prime}(D)$   
as $\frac{1}{2}-\frac{1}{D} > 0$}\\
& \leq & \frac{C^{\prime \prime}}{\eta} \ \ \ \mbox{ as $\eta \geq   
\sqrt{ D \log \log n}$} 
\end{eqnarray*}
\end{proof}

\subsection{Proofs for Section~\ref{sec:Optimality-of-spatial-HC-BJ}}

\begin{proof}[Proof of Theorem~\ref{thm: Spatial HC is optimal}](i)
We first prove optimality for $sHC_{n}$ and then derive the
conclusion for the other statistics from this result.

We will show that if $r>\rho^{*}(\alpha,\beta)$, then
$sHC_{n}=\Omega_{p}(n^{\xi})$ for some $\xi>0$. Then the claim
about $sHC_{n}$ follows with the result about the null distribution
given in Theorem~\ref{thm:HC-BJ-Under-the-null}.

Let $\ell^{*}$ be the level that corresponds to the true length of
the signal, i.e. $\ell^{*}$ satisfies $2^{\ell^{*}-1}<n^{\alp}\le2^{\ell^{*}}$.
Note that  $0\le\alpha<1$ implies $\lmax-\ell^{*}=\Theta(\log n)$. 
Further, Lemma~\ref{lem: Short-range-dependent} shows that $n_{\ell^{*}}:=\#\Iapp (\ell^{*})$
satisfies 
\be \label{boundn}
\frac{1}{2} n^{1-\alpha}\ \leq \ n 2^{-\ell^*} \ \leq \ n_{\ell^{*}}\ \leq \ 144 n 2^{-\ell^*} \log_2n \  
\leq \ 144 n^{1-\alpha} \log_2 n.
\ee
Below we will consider the 
two disjoint situations $r/(1-\alpha)<\frac{1}{4}$ and $r/(1-\alpha)\ge\frac{1}{4}$. We define
$t^{*}$ such that 

\[
\bar{\Phi}^{-1}(t^{*})=\begin{cases}
2\sqrt{2r\log n} & r/(1-\alpha)<\frac{1}{4}\\
\sqrt{2\log n_{\ell^{*}}-5\log\log n_{\ell^{*}}} & r/(1-\alpha)\ge\frac{1}{4}
\end{cases}.
\]
Solving for $t^{*}$, we have 

\[
t^{*}=\begin{cases}
L_{n}n^{-4r} & r/(1-\alpha)<\frac{1}{4}\\
\frac{\log^{2}n_{\ell^{*}}}{4\sqrt{\pi}n_{\ell^{*}}}(1+o(1)) & r/(1-\alpha)\ge\frac{1}{4}
\end{cases}.
\]
By (\ref{represent})
\begin{eqnarray}
sHC_{n} & \ge & \sqrt{\frac{n}{2^{\ell^{*}}n_{\ell^{*}}}} HC_{n_{\ell^{*}}}(\ell^{*}) \nn \\
&  = & \sqrt{\frac{n}{2^{\ell^{*}}n_{\ell^{*}}}}\sup_{t\in [p_{(1)},\, p_{(n_{\ell^{*}}/2)}]}
\frac{\sum{}_{I\in
\Iapp(\ell^{*})}\left({\bf 1}(\X (I)\ge\bar{\Phi}^{-1}(t))-t\right)}{
\sqrt{n_{\ell^{*}}t(1-t)}}  \label{revlabel}
\end{eqnarray}

On the event $\{ p_{(n_{\ell^{*}}/2)} <t^{*}\}$ we set $t:=p_{(n_{\ell^{*}}/2)}$ to see that the sup is not
smaller than
$$
\frac{\sum_{I\in
\Iapp(\ell^{*})}\Bigl({\bf 1}\left(\bar{\Phi}(\X (I))\le t \right)-t \Bigr)}{
\sqrt{n_{\ell^{*}} \,t \left(1-t \right)}} \ =\ 
\frac{ \frac{ n_{\ell^{*}}}{2} - n_{\ell^{*}} \,p_{(n_{\ell^{*}}/2)} }{
\sqrt{n_{\ell^{*}} \,p_{(n_{\ell^{*}}/2)} \left(1-p_{(n_{\ell^{*}}/2)} \right)}} \ \geq \
\frac{n_{\ell^{*}} (\frac{1}{2} -t^{*})}{\sqrt{n_{\ell^{*}} t^{*} (1-t^{*})}}
$$
(\ref{revlabel}), (\ref{boundn}) and $t^{*} =o(1)$ show that $sHC_{n} \geq n^{\frac{1-\alp}{2}}$  
for $n$ large enough.

Now we consider the event $\{ p_{(n_{\ell^{*}}/2)} \geq t^{*}\}$. We will show below
\begin{equation}
\P_{H_{1}}\Bl(p_{(1)}>\frac{\log^{3/2}n_{\ell^{*}}}{n_{\ell^{*}}}\Br)\rightarrow 0.
\label{eq: p_1 not too small}
\end{equation}
But if $p_{(1)}\le\frac{\log^{3/2}n_{\ell^{*}}}{n_{\ell^{*}}}$,
then we have $t^{*}\ge p_{(1)}$ for $n$ large enough by (\ref{boundn}). On the event
$\{p_{(1)} \leq t^{*} \leq  p_{(n_{\ell^{*}}/2)} \}$ we obtain from (\ref{revlabel})
$$
sHC_{n} \ \geq \ \sqrt{\frac{n}{2^{\ell^{*}}n_{\ell^{*}}}}\, \frac{\sum{}_{I\in\Iapp(\ell^{*})}
\left({\bf 1}(\X (I)\ge\bar{\Phi}^{-1}(t^{*}))-t^{*}\right)}{\sqrt{n_{\ell^{*}}t^{*}(1-t^{*})}}
 \ =: \ T_n(\ell^*)
$$
We will show that $\E\,T_n(\ell^*) =\Omega(n^{\xi})$ for some $\xi>0$ and
$\sqrt{\V\,T_n(\ell^*)} =o(\E\,T_n(\ell^*))$. Then Chebychev's inequality will yield the
desired conclusion
\begin{equation}
sHC_{n}\ \geq \ T_n(\ell^*) \ =\ \Omega_{p}(n^{\xi}).\label{eq: hc_alter_1}
\end{equation}

Recall the notation $\X(I):=\sum_{i \in I}X_i /\sqrt{|I|}$, so $\X (I) \sim \mathcal{N}(\E \X (I),1)$.
Denote $\mu':=\sqrt{2r\log n}(1-\frac{1}{3\sqrt{\lmax-\ell^{*}+4}})$.
By the  construction of $\Iapp(\ell^{*})$ (see also Proposition~\ref{properties}(ii)), 
there are at least
$n^{1-\alpha-\beta}$ intervals $I \in \Iapp(\ell^{*})$  satisfying 
\be \label{meanlower}
\E(\X(I))\ge\mu'=\sqrt{2r\log n}-O(1).
\ee
Situation 1: If $r/(1-\alpha)<\frac{1}{4}$, then we have:

\begin{eqnarray*}
\E\,T_n(\ell^*) & \ge & \sqrt{\frac{n}{2^{\ell^{*}}
 n_{\ell^{*}}}}\frac{1}{\sqrt{n_{\ell^{*}}t^{*}(1-t^{*})}}
\,n^{1-\alpha-\beta}\Bl(\P \bigl(\mathcal{N}(\mu',1)\ge\bar{\Phi}^{-1}(t^{*})\bigr)-t^{*}\Br)\\
& = & L_n n^{\frac{1-\alpha}{2}-\beta+2r} \Bl(\P\bigl(\mathcal{N}(0,1) \ge 2 \sqrt{2r \log n}
  -\mu'\bigr) -t^*\Br)\\
 & \ge & L_{n}n^{\frac{1-\alpha}{2}-\beta+r}
\end{eqnarray*}

by (\ref{boundn}) and Mill's ratio. 
$\rho^{*}(\alpha,\beta)<r<\frac{1-\alpha}{4}$ implies $\frac{1-\alpha}{2}-\beta+r>0$,
so we can take $0<\xi<\frac{1-\alpha}{2}-\beta+r$ to conclude $\E\,T_n(\ell^*)=\Omega(n^{\xi})$.

In order to compute the variance of $T_n(\ell^*)$ note that by Lemma~\ref{lem: Short-range-dependent} 
the intervals in $\Iapp (\ell^*)$ can be grouped into $i_{\mathrm{max}} \leq 144 \log_2 n$ groups
$\mathcal{J}_{i}(\ell^{*})$, $i=1,\ldots,i_{\mathrm{max}}$,
each of which contains not more than $\# \Iapp (\ell^*) =n_{\ell^{*}}$ disjoint intervals.
Thus within each group $\mathcal{J}_{i}(\ell^{*})$ the $\X(I)$ are independent and therefore

\begin{align*}
 \V\,  \Bl(\sum_{I\in\mathcal{J}_{i}(\ell^{*})} & {\bf 1}\bigl(\X(I)\ge\bar{\Phi}^{-1}(t^{*})\bigr)\Br)
  =  \sum_{I\in\mathcal{J}_{i}(\ell^{*})}\P\bigl(\X(I)\ge\bar{\Phi}^{-1}(t^{*})\bigr)
\bigl(1-P(\X(I)\ge\bar{\Phi}^{-1}(t^{*}))\bigr)\\
 & \le  \sum_{I\in\Iapp(\ell^{*}),E\X(I)=0}\P(\X(I)\ge\bar{\Phi}^{-1}(t^{*}))
\bigl(1-\P(\X(I)\ge\bar{\Phi}^{-1}(t^{*}))\bigr)\\
 &   \hspace{0.5cm} +\sum_{I\in\mathcal{J}_{i}(\ell^{*}),E\X(I)>0}\P(\X(I)\ge\bar{\Phi}^{-1}(t^{*}))\\
 & \le  L_{n}(n^{1-\alpha-4r}+n^{1-\alpha-\beta-r})
\end{align*}
by (\ref{boundn}) and since the number of $I \in \mathcal{J}_{i}(\ell^{*})$ that intersect with one 
of the $m=n^{1-\alpha -\beta}$ intervals that have an elevated mean can not be larger than $2m$,
and an overlap results in $\E \X(I) \leq \sqrt{2r \log n}$.
Applying Cauchy-Schwartz to the covariances
between the $i_{\mathrm{max}} \leq 144 \log_2 n$ groups gives 
$$
\V \, T_n(\ell^*)\ \leq \ \frac{n}{2^{\ell^{*}} n_{\ell^{*}}^2 t^*(1-t^*)}
\ i_{\mathrm{max}}^2 L_{n}(n^{1-\alpha-4r}+n^{1-\alpha-\beta-r}) \ \leq \ L_n (1+n^{3r-\beta})
$$
by (\ref{boundn}). 
Since $\rho^{*}(\alpha,\beta)<r<\frac{1-\alpha}{4}$
implies $\max(0,\frac{3r-\beta}{2})<\frac{1-\alpha}{2}-\beta+r$,
we conclude $\sqrt{\V \, T_n(\ell^*)} =o(\E \,T_n(\ell^*))$, and thus 
Chebychev's inequality gives (\ref{eq: hc_alter_1}).
\smallskip

Situation 2: When $r \ge \frac{1-\alpha}{4}$, by a very similar calculation as above
we obtain
$\E \,T_n(\ell^*) \ge L_{n}n^{1-\alpha-\beta}n^{-(\sqrt{1-\alpha}-\sqrt{r})^{2}}$ and
$\sqrt{\V \, T_n(\ell^*)} =o(\E \,T_n(\ell^*))$.
 Since we assume $r>\max(\frac{1-\alpha}{4},\allowbreak \rho^{*}(\alpha,\beta)$),
we can find $\xi$ with  $0<\xi<1-\alpha-\beta-(\sqrt{1-\alpha}-\sqrt{r})^{2}$,
and hence (\ref{eq: hc_alter_1}) also follows in this situation.

Thus we have shown that $r>\rho^{*}(\alpha,\beta)$ and 
$p_{(1)}\le\frac{\log^{3/2}n_{\ell^{*}}}{n_{\ell^{*}}}$ imply (\ref{eq: hc_alter_1})
for some $\xi >0$, and the proof for $sHC_{n}$ will be complete once we show
(\ref{eq: p_1 not too small}):

\begin{eqnarray*}
\P_{H_{1}}\Bl(p_{(1)}>\frac{\log^{3/2}n_{\ell^{*}}}{n_{\ell^{*}}}\Br)
 & = & \P_{H_{1}}\Bl(\max_{I\in\Iapp (\ell^{*})}\X(I)<\bar{\Phi}^{-1}(
\frac{\log^{3/2}n_{\ell^{*}}}{n_{\ell^{*}}})\Br)\\
 & = & \P_{H_{1}}\Bl(\max_{i=1}^{i_{\mathrm{max}}}\max_{I\in\mathcal{J}_{i}(\ell^{*})}\X(I)
<\bar{\Phi}^{-1}(\frac{\log^{3/2}n_{\ell^{*}}}{n_{\ell^{*}}})\Br)\\
 & \le & \min_{i=1}^{i_{\mathrm{max}}}\P_{H_{1}}\Bl(\max_{I\in\mathcal{J}_{i}(\ell^{*})}\X(I)
<\bar{\Phi}^{-1}(\frac{\log^{3/2}n_{\ell^{*}}}{n_{\ell^{*}}})\Br)\\
 & \le & \min_{i=1}^{i_{\mathrm{max}}}\P_{H_{0}}\Bl(\max_{I\in\mathcal{J}_{i}(\ell^{*})}\X(I)
<\bar{\Phi}^{-1}(\frac{\log^{3/2}n_{\ell^{*}}}{n_{\ell^{*}}})\Br)\\
 & \le & \Bl(1-\frac{\log^{3/2}n_{\ell^{*}}}{n_{\ell^{*}}}\Br)^{n2^{-\ell^{*}}-2}
\end{eqnarray*}
since $\# \mathcal{J}_{i}(\ell^{*}) \geq n2^{-\ell^{*}}-2$ by Lemma~\ref{lem: Short-range-dependent}.
Now (\ref{eq: p_1 not too small}) follows since
$\frac{n2^{-\ell^{*}}\log^{3/2}n_{\ell^{*}}}{n_{\ell^{*}}}\ge\frac{\log^{3/2}
n_{\ell^{*}}}{144 \log_2 n}\rightarrow+\infty$ by (\ref{boundn}), completing the proof for $sHC_{n}$.

As for $sS_{n}^{}(s)$, Lemma 7.2 in \cite{jager2007goodness} shows that 
$K_{s}(u,v) 1(v<u \leq \frac{1}{2}) \le K_{2}(u,v) 1(v<u\leq \frac{1}{2})$ for all $s \in [-1,2]$. 
Thus $S_{n}^{+}(s)\le S_{n}^{+}(2) = \frac{1}{2} (HC_n^+)^2$ and therefore
$sS_{n}^{}(s)\le \frac{1}{2} \left(sHC_n^{} \right)^2$. 
Hence it follows from Theorem~\ref{thm:HC-BJ-Under-the-null} that under the null distribution
\begin{equation}
sS_{n}^{}(s)=O_{p}(\log^{4}n)\label{eq:spatial-goodness-of-fit-null}
\end{equation}
for all $-1\le s\le2$. (That theorem also provides a better bound for the special case $s=1$.)

Now we examine the performance of $sS_{n}^{}(s)$ when $r>\rho^{*}(\alpha,\beta)$.
As in \cite{donoho2004higher}, we need to consider two cases: $\rho^{*}(\alpha,\beta)<r<\beta/3$
and $r>(\sqrt{1-\alpha}-\sqrt{1-\alpha-\beta})^{2}$. These two cases
overlap and together cover the full region $r>\rho^{*}(\alpha,\beta)$.

In the first case where $0<\mbox{\ensuremath{\rho}}^{*}(\alpha,\beta)<r<\beta/3$,
we must have $\beta<\frac{3}{4}(1-\alpha)$ and hence $r<\frac{1-\alpha}{4}$, and so
we can choose a positive $r_{0}<r<\frac{1-\alpha}{4}$.
Let $\ell^{*}$ be the level that corresponds to the true length of
the signal. Define

\[
HC_{n_{\ell^{*}},r,r_{0}}(\ell^{*})\ :=
\sup_{n^{-4r}\le p_{(i)}\le n^{-4r_{0}}}
\sqrt{n_{\ell^{*}}}\frac{\frac{i}{n_{\ell^{*}}}-p_{(i)}}{\sqrt{p_{(i)}(1-p_{(i)})}}
{\bf 1}\Bl(p_{(i)} < \frac{i}{n_{\ell^{*}}}\Br)
\]
where the p-values pertain to intervals in $\Iapp(\ell^{*})$. We
need the following lemma which is proved in the \hyperref[appn]{Appendix}:

\begin{lemma}
\label{lem: goodness_of_fit_converge_to_one}
Let $p_{(i)}$ be the ordered p-values for intervals in $\Iapp(\ell^{*})$.
Then $0<\mbox{\ensuremath{\rho}}^{*}(\alpha,\beta)<r<\beta/3$ implies
$\ \sup_{n^{-4r}<p_{(i)}<n^{-4r_{0}}}
|\frac{i}{n_{\ell^{*}}p_{(i)}}-1|\stackrel{p}{\rightarrow}0$.
\end{lemma}

Using the above lemma and Lemma 7.2 in \cite{jager2007goodness}
we have 
\[
\sup_{n^{-4r}<p_{(i)}<n^{-4r_{0}}}n_{\ell^{*}}K_{s}\Bl(\frac{i}{n_{\ell^{*}}},
p_{(i)}\Br)
{\bf 1}\Bl(p_{(i)} < \frac{i}{n_{\ell^{*}}}\Br) \ge
\frac{1}{2}(HC_{n_{\ell^{*}},r,r_{0}}(\ell^{*}))^{2}(1+o_{p}(1))
\]
Thus
\begin{eqnarray}
sS_{n}^{}(s) & \ge & \frac{n}{2^{\ell^{*}}n_{\ell^{*}}}S_{n_{\ell^{*}}}^{+}(s,\ell^{*})\nonumber \\
 & \ge & \frac{n}{2^{\ell^{*}}n^{\ell^{*}}}\sup_{n^{-4r}<p_{(i)}<n^{-4r_{0}}}
n_{\ell^{*}}
K_{s}\Bl(\frac{i}{n_{\ell^{*}}},p_{(i)}\Br) {\bf 1}\Bl(p_{(i)} < \frac{i}{n_{\ell^{*}}}\Br) \nonumber \\
 & \ge & \frac{1}{2}\Bl(\sqrt{\frac{n}{2^{\ell^{*}}n^{\ell^{*}}}}HC_{n_{\ell^{*}},r,r_{0}}
(\ell^{*})\Br)^{2}(1+o_{p}(1))=\Omega_{p}(n^{\xi})\label{eq: spatial alternative - 1}
\end{eqnarray}
for some $\xi>0$ by the above proof about $sHC_{n}$ that localized the analysis to 
$t^*=L_n n^{-4r}$.

For the second case, if $r>(\sqrt{1-\alpha}-\sqrt{1-\alpha-\beta})^{2}$ and $r<1-\alpha$,
then $(r+\beta)/2\sqrt{r}<\sqrt{1-\alpha}$. So we can pick $q \in (0,1)$
such that $ $$\max((r+\beta)/2\sqrt{r},\sqrt{r})<\sqrt{q}<\sqrt{1-\alpha}$.
As noted above, there are at least $n^{1-\alpha-\beta}$ intervals $I \in \Iapp(\ell^{*})$
satisfying (\ref{meanlower}).
By Lemma \ref{lem: Short-range-dependent} the intervals in $\Iapp(\ell^{*})$
can be grouped into at most $144\log_2 n$ groups such
that each group consists of disjoint
intervals. By the pigeonhole principle, at least one group contains
more than $\frac{n^{1-\alpha-\beta}}{144 \log_2 n}$ intervals satisfying (\ref{meanlower}).
Since the $\X(I)$ in that group are independent we have

\be \label{binlower}
\sum_{I\in\Iapp(\ell^{*})}\boldsymbol{1}(\X(I)\geq \sqrt{2q\log n})\ \stackrel{d}{\ge} \
{\rm Bin}\Bl(\frac{n^{1-\alpha-\beta}}{144\log_2 n},L_{n}n^{-(\sqrt{q}-\sqrt{r})^{2}}\Br)
\ =\ \Omega_{p}(n^{\xi})
\ee
by Chebychev's inequality, since $\xi:=1-\alpha-\beta-(\sqrt{q}-\sqrt{r})^2 >1-\alpha-q >0$.
Setting $t:=\sqrt{2q \log n}$ we get
$$
\frac{\sum_{I\in\Iapp(\ell^{*})} \boldsymbol{1}
(\X(I)\geq t)}{n_{\ell^{*}}\bar{\Phi}(t)} \stackrel{p}{\rightarrow}\infty
$$
since $n_{\ell^{*}} \bar{\Phi}(t) =L_n n^{1-\alpha -q}$ by (\ref{boundn}) and $\xi>1-\alpha-q$.
 
Together with Lemma~7.2 in \cite{jager2007goodness} and (\ref{represent}) we obtain
\begin{eqnarray}
S_{n_{\ell^{*}}}^{+}(s,\ell^{*}) & \geq & n_{\ell^{*}} 
K_s\Bl( \frac{\sum_{I\in\Iapp(\ell^{*})} \boldsymbol{1} (\X(I)\geq t)}{n_{\ell^{*}}},
\bar{\Phi}(t) \Br) {\bf 1} \Bl( \bar{\Phi}(t) < 
\frac{\sum_{I\in\Iapp(\ell^{*})} \boldsymbol{1} (\X(I)\geq t)}{n_{\ell^{*}}}\Br) \nn \\
& \geq & L_{n}(1+o_p(1))\sum_{I\in\Iapp(\ell^{*})}\boldsymbol{1}(\X(I)\geq t)\nn.
\end{eqnarray}
It follows with (\ref{binlower}) and (\ref{boundn}) that
\begin{equation}
sS_{n}^{}(s)\ \ge \ \frac{n}{2^{\ell^*}n_{\ell^{*}}}S_{n_{\ell^{*}}}^{+}(s,\ell^{*})\ =\
\Omega_{p}(n^{\xi}). \label{eq: spatial alternative -2}
\end{equation}
From equations \eqref{eq:spatial-goodness-of-fit-null}, \eqref{eq: spatial alternative - 1}
 and \eqref{eq: spatial alternative -2}
it follows that for all $-1\le s\le2$, $sS_{n}^{}(s)$ has asymptotic power 1  under the alternative
$r>\rho^{*}(\alpha,\beta)$.
\end{proof}

\begin{proof}[Proof of Theorem~\ref{thm: Spatial HC is optimal}](ii)
The null case was discussed in Theorem \ref{thm:HC-BJ-Under-the-null} which showed
$sHC_{n}=O_{p}(\log^{2}n)  
$. As in the proof of part (i), we only
need to show that $sHC_{n}=\Omega_{p}(n^{\xi})$ for some $\xi>0$
when $r>\rho^{*}(\alpha,\beta)$. Again, let $\ell^{*}$ be the level
that corresponds to the true length of the signal, i.e. $\ell^{*}$
satisfies $2^{\ell^{*}-1}< n^{\alp} \le2^{\ell^{*}}$. $0\le\alpha<1$ implies
$\lmax-\ell^{*}=\Theta(\log n)$. By the construction of $\Iapp(\ell^*)$
(see also Proposition~\ref{properties}(ii)) there are at least $n^{1-\alpha-\beta}$
intervals $I \in \Iapp(\ell^{*})$  satisfying
$\E(\X(I))\ge\ n^r (1-\frac{1}{3\sqrt{\lmax-\ell^{*}+4}}) \geq n^r/2$. Hence

\bean
\E \left(\sum_{I \in \Iapp(\ell^*)} \Bl( {\bf 1} (\X(I) \geq \bar{\Phi}^{-1}(\frac{1}{4}))
  -\frac{1}{4}\Br)\right)
& \ge & n^{1-\alpha-\beta} \Bl(\bar{\Phi}(\bar{\Phi}^{-1}(\frac{1}{4})- n^r/2)-\frac{1}{4}\Br)\nn\\
& \geq & n^{1-\alpha-\beta} n^r/8  \label{densestar}
\eean
since $\bar{\Phi}' \leq -\frac{1}{4}$ on $\Bl( \bar{\Phi}^{-1}(\frac{1}{4})-\frac{1}{2},
\bar{\Phi}^{-1}(\frac{1}{4})\Br)$ and we may w.l.o.g. assume that $n^r <1$ since $\rho^{*}<0$.

By Lemma~\ref{lem: Short-range-dependent}, the intervals in $\Iapp (\ell^*)$ can be grouped into
at most $144 \log_2 n$
groups, each of which contains not more than $\# \Iapp (\ell^*) =n_{\ell^{*}}$ disjoint intervals.
Thus within each group the $\X(I)$ are independent, and applying Cauchy-Schwartz to the covariances
between groups gives
\bea
\V \Bl(\sum_{I \in \Iapp(\ell^*)} \Bl( {\bf 1} (\X(I) \geq \bar{\Phi}^{-1}(\frac{1}{4}))-
   \frac{1}{4}\Br)\Br)
& \leq & (144 \log_2 n)^2 n_{\ell^{*}} \V\Bl({\bf 1} (\X(I) \geq \bar{\Phi}^{-1}(\frac{1}{4}))\Br)\\
& \leq & L_n n_{\ell^{*}}
\eea
Together with (\ref{densestar}) and (\ref{boundn}) this shows that
$$
T_n(\ell^*)\ :=\ \sqrt{\frac{n}{2^{\ell^{*}}n_{\ell^{*}}}}
 \frac{\sum_{I \in \Iapp(\ell^*)} \Bl( {\bf 1} (\X(I) \geq \bar{\Phi}^{-1}(\frac{1}{4}))-\frac{1}{4}\Br)}{
  \sqrt{n_{\ell^{*}}\frac{1}{4} (1-\frac{1}{4})}}
$$
satisfies $\E\, T_n(\ell^*) \geq L_n n^{\frac{1-\alp}{2} -\beta +r}$ and $\V\, T_n(\ell^*)
\leq L_n$, hence Chebychev and $r>\beta -\frac{1-\alpha}{2}$ yield
\be \label{dense2star}
T_n(\ell^*) \ =\ \Omega_p (n^{\xi})
\ee
for some $\xi>0$.

Now we partition the sample space into three events:
\begin{eqnarray*}
sHC_{n} & \ge & HC_{n_{\ell^{*}}}(\ell^{*}) \Bl( {\bf 1} \Bl(p_{(n_{\ell^{*}}/2)}<\frac{1}{4}\Br)
  +  {\bf 1} \Bl(p_{(1)} \leq \frac{1}{4} \leq p_{(n_{\ell^{*}}/2)}\Br)
  +  {\bf 1} \Bl( \frac{1}{4} < p_{(1)}\Br) \Br)\\
 & \ge & \sqrt{n_{\ell^*}} \frac{\frac{1}{2}-\frac{1}{4}}{\sqrt{\frac{1}{4}(1-\frac{1}{4})}}
  \,{\bf 1} \Bl(p_{(n_{\ell^{*}}/2)}<\frac{1}{4}\Br)
   +T_n(\ell^*) \,{\bf 1} \Bl(p_{(1)} \leq \frac{1}{4} \leq p_{(n_{\ell^{*}}/2)}\Br)\\
   & & +HC_{n_{\ell^{*}}}(\ell^{*}) \,{\bf 1} \Bl( \frac{1}{4} < p_{(1)}\Br)\\
& =& \Omega_p \Bl(n^{\min (\frac{1-\alpha}{2},\xi)}\Br)
\eea
by (\ref{dense2star}), (\ref{boundn}) and ${\rm P } (\frac{1}{4} <p_{(1)}) \leq (\frac{3}{4})^{n_{\ell^{*}}}
\rightarrow 0$.
\end{proof}

\subsection{Proofs for Section~\ref{sec:Comparison-with-other-statistic}}

\begin{proof}[Proof of Theorem \ref{thm:detection-boundary-penalized-scan}]
In the dense case let $r:=\beta - \frac{1-\alpha}{2}+\epsilon$ for
some $\eps >0$. Then $\sum_{i=1}^n X_i/\sqrt{n}$ is normal with variance
one and mean $n^{1-\beta} \mu /\sqrt{n}= n^{1-\beta +r -\frac{1+\alp}{2}}
=n^{\eps}$. Thus $P_n \geq \sum_{i=1}^n X_i/\sqrt{n} -\sqrt{2} \stackrel{p}{
\rightarrow} \infty$. Hence $P_n$ has asymptotic power one since $P_n=O_p(1)$
under $H_0$.

In the sparse case,
if $r>\rho_{\mathrm{pen}}^{*}(\alpha,\beta)$ then we can pick a constant
$\epsilon>0$ depending only on $(r,\alpha,\beta)$ such that 
$1-\alpha-\beta>((1+\epsilon)\sqrt{1-\alpha}-\sqrt{r})^{2}$.
For a block $I_g=:(j,j+n^{\alp}]$ in (\ref{eq: Def-Block-mixture-problem}) write 
$Z_g:=\sum_{i=j+1}^{j+n^{\alp}} X_i/\sqrt{n^{\alp}}$.
In order to show that $P_n$ has asymptotic power one, it is enough to show that
\be  \label{thm9st}
\P_{\mu (n)} \Bl( \max_{g=1,\ldots,m} Z_g\ >\ (1+\eps)\sqrt{2 \log \frac{n}{n^{\alp}}}
\Br) \rightarrow 1
\ee
because $P_n =O_p(1)$ under $H_0$ while $\eps \sqrt{2 \log \frac{n}{n^{\alp}}}
=\eps \sqrt{2 (1-\alp)\log n} \rightarrow \infty$.

Note that the $Z_g$ are independent normal with mean $\sqrt{n^{\alp}} \mu
=\sqrt{2r \log n}$ and variance one. Therefore
\begin{align*}
p_n & := \P\Bl(Z_1 > (1+\eps)\sqrt{2 \log \frac{n}{n^{\alp}}} \Br) \\
& = 1 -\Phi \Bl((1+\epsilon) \sqrt{2 (1-\alp)\log n} -\sqrt{2r \log n} \Br) \\
& \geq L_n n^{-(\sqrt{r} -(1+\eps) \sqrt{1-\alp})^2}
\end{align*}
by Mill's ratio. Hence the probability in (\ref{thm9st}) equals
$$
1-(1-p_n)^m \geq 1-\exp (-mp_n)
$$
and $mp_n \geq L_n n^{1-\alp -\beta -(\sqrt{r} -(1+\eps) \sqrt{1-\alp})^2}
\rightarrow \infty$.

The claim for $P_n^{\mathrm{app}}$ obtains in the same way, by taking account of the
approximation error incurred by using the approximating set, see
\cite[Theorem~2]{rivera2013optimal} and \cite[Theorem~11]{KouThesis}. 

Proceeding as in the proof of Theorem~1.4 in \cite{donoho2004higher},
it can be shown that $P_n$ and $P_n^{{\mathrm{app}}}$ are powerless if $r<\rho_{\mathrm{pen}}^{*}(\alpha,\beta)$. 
\end{proof}

\subsection{Proofs for Section~\ref{multivariate}}

The following lemma is the bivariate analogue of Lemma~\ref{lem: Short-range-dependent}:

\begin{lemma} \label{lemmamult}
The rectangles in $\Iapp^{(2)} (\ell)$ can be grouped into at most $12 \eps_{\ell}^{-4}
(\ell+1) \leq 8 \cdot 6^5 (\log_2 n)^2 (\ell+1)$ groups such that each group
consists of at least $\frac{9}{16} \frac{n^2}{2^{\ell}}$ and at most $2 \frac{n^2}{2^{\ell}}$
disjoint rectangles. Hence $\# \Iapp^{(2)} (\ell) \leq 16 \cdot 6^5 (\log_2 n)^2 n^2 
\frac{\ell +1}{2^{\ell}}$.
\end{lemma}

\begin{proof}[Proof of Lemma~\ref{lemmamult}] We will use the following refinement of
Lemma~\ref{lem: Short-range-dependent} for the univariate setting:

\begin{claim} \label{claim1}
 The intervals in $\Iapp(\ell)$ that have a given length $L$ (which
hence is a multiple of $d_{\ell}$) can be grouped into $\frac{L}{d_{\ell}}$ groups such that
each group consists of either $\lfloor \frac{n}{L} \rfloor$ or $\lfloor \frac{n}{L}\rfloor -1$
disjoint intervals.
\end{claim}

To see this, set $I_j:=(jd_{\ell},jd_{\ell} +L]$ for $j=0,\ldots,\frac{L}{d_{\ell}}-1$,
and consider all possible shifts of $I_j$ by multiples of $L$:
$$
C(j)\ :=\ \Bigl\{ kL+I_j, \ k=0,\ldots,\lfloor \frac{n-jd_{\ell}}{L}\rfloor -1 \Bigr\}
$$
One readily checks that $\bigcup_{j=0}^{\frac{L}{d_{\ell}}-1} C(j)$ equals the collection of
all intervals in $\Iapp(\ell)$ that have length $L$. Further, each $C(j)$ consists of
$\lfloor \frac{n}{L} \rfloor$ or $\lfloor \frac{n}{L}\rfloor -1$ intervals that are
disjoint, proving Claim~\ref{claim1}.
\smallskip

Now we consider the rectangles in $\Iapp^{(2)}(\ell)$
that have given sidelengths $L_1$ and $L_2$:
\begin{align*}
C(\ell,L_1,L_2) &:=\Bigl\{ R=I_1 \times I_2 \in \Iapp^{(2)}(\ell): |I_1|=L_1, |I_2|=L_2 \Bigr\}\\
& = \Bigl\{ R=I_1 \times I_2 \in \Iapp(\ell_1,\eps_{\ell}) \times \Iapp(\ell_2,\eps_{\ell}):
|I_1|=L_1, |I_2|=L_2 \Bigr\}
\end{align*}
where $\ell_i = \lceil \log_2 L_i \rceil$. 

\begin{claim} \label{claim2}
The rectangles in $C(\ell,L_1,L_2)$ can be grouped into at most 
$4 \eps_{\ell}^{-2} \leq 4 \cdot 6^2 \log_2 n^2$ groups such that each group
consists of at least $(\lfloor \frac{n}{L_1} \rfloor -1) (\lfloor \frac{n}{L_2} \rfloor -1)
\geq \frac{9}{16} \frac{n^2}{2^{\ell}}$  and at most
$\lfloor \frac{n}{L_1} \rfloor  \lfloor \frac{n}{L_2} \rfloor \leq 2 \frac{n^2}{2^{\ell}}$
disjoint rectangles.
\end{claim}

In order to prove Claim~\ref{claim2}, note that Claim~\ref{claim1} implies that the rectangles in 
$C(\ell,L_1,L_2)$ can be grouped into $\frac{L_1}{d_{\ell_1}} \times \frac{L_2}{d_{\ell_2}}
\leq \frac{4}{\eps_{\ell}^2} \leq 4 \cdot 6^2 \log_2 n^2$ groups such that each group contains
between $(\lfloor \frac{n}{L_1}-1 \rfloor)(\lfloor \frac{n}{L_2}-1 \rfloor)$ and
$\lfloor \frac{n}{L_1} \rfloor \lfloor \frac{n}{L_2} \rfloor$ rectangles that are
disjoint (since the Cartesian product of two collections of disjoint intervals yields
a collection of disjoint rectangles). Since the area of the rectangles satisfies
$L_1L_2 \in (2^{\ell-1},2^{\ell}]$ we get $\lfloor \frac{n}{L_1} \rfloor \lfloor 
\frac{n}{L_2} \rfloor \leq 2 \frac{n^2}{2^{\ell}}$. Finally, $L_i \leq n/8$ implies
$(\lfloor \frac{n}{L_1}-1 \rfloor)(\lfloor \frac{n}{L_2}-1 \rfloor) \geq (\frac{3}{4}
\frac{n}{L_1})(\frac{3}{4} \frac{n}{L_2}) \geq \frac{9}{16} \frac{n^2}{2^{\ell}}$,
establishing Claim~\ref{claim2}.
\smallskip

The lemma now obtains as follows:
Clearly, $\Iapp^{(2)} (\ell) =\bigcup_{\{\mbox{all possible $L_1,L_2$}\}} C(\ell,L_1,L_2)$.
Since the level $\ell_1$ of $L_1$ must satisfy $\ell_1 \leq \ell$ and each $\Iapp(\tilde{\ell},
\eps_{\ell})$ admits at most $\lceil 2^{\tilde{\ell}-1}/d_{\tilde{\ell}} \rceil \leq
\lceil 1/\eps_{\ell} \rceil$ different interval lengths, there are at most  $\lceil
1/\eps_{\ell} \rceil (\ell +1)$ different choices for $L_1$. The constraint $\ell \leq
\ell_1+\ell_2 \leq \ell +1$ from Proposition~\ref{properties}
implies that given $L_1$, the level $\ell_2$ of $L_2$ must
be either $\ell -\ell_1$ or $\ell -\ell_1+1$, hence there are at most $\lceil 2/\eps_{\ell}
\rceil$ different choices for $L_2$. So there are at most $\frac{3}{\eps_{\ell}^2}(\ell +1)
\leq 3 \cdot 6^2  (\log_2 n^2)(\ell +1)$ different choices for $(L_1,L_2)$. Lemma~\ref{lemmamult}
now follows with Claim~\ref{claim2}. We note that the statement of the lemma can be sharpened somewhat
as the factor $\frac{9}{16}$ is due to large rectangles which allow a better bound on $12 \eps_{\ell}^{-4}
(\ell+1)$. 
\end{proof}

\begin{proof}[Proof of Theorem~\ref{thm:multivariate-null}] The proof follows that of 
Theorem~\ref{thm:HC-BJ-Under-the-null} using the inequalities from Lemma~\ref{lemmamult} 
in place of Lemma~\ref{lem: Short-range-dependent}. That is, for a fixed $\ell$ we now have
$n_{\ell} \leq 16 \cdot 6^5 (\log_2 n)^2 n^2$, $\frac{9}{16} \frac{n^2}{2^{\ell}} \leq \# G_i \leq
n^2$, $i_{\mathrm{max}} \leq 16 \cdot 6^5 (\log_2 n)^3$ and $\lmax +1 \leq 2 \log_2 n$. 
As for $sBJ_{n}^{(2)}$, the two
additional factors of $\log_2 n$ in $i_{\mathrm{max}}$ and the factor $\frac{9}{16}$ in the lower bound
for $\# G_i$ necessitate to replace the condition $c>3$ by $c>(3+2)\frac{16}{9}$ in order
to obtain the desired convergence to 0. This bound on $c$ can be improved somewhat by refining
the bounds in Lemma~\ref{lemmamult} as explained at the end of its proof.
Concerning $sHC_{n}^{(2)}$, the convergence rate needs to account for the
two additional factors of $\log_2 n$ in $i_{\mathrm{max}}$. 
\end{proof}

\begin{proof}[Proof of Theorem~\ref{optimality2}] The proof of the lower bound is analogous to that of 
Theorem~\ref{thm: Block-mixture-detection-boundary}
by considering the submodel obtained by partitioning the $n \times n$ grid into 
$n'=n^{2-2\alp}$ blocks of size $|I|=n^{2\alp}$. The claim about $sHC_{n}^{(2)}$
and $sBJ_{n}^{(2)}$
follows as in Theorem~\ref{thm: Spatial HC is optimal} by using $n^2$ in place of $n$. 
\end{proof}

\subsection{Proofs for Section~\ref{graphs}}

\begin{proof}[Proof of Proposition~\ref{PropC}] There are at most $\frac{n}{d_{\ell}} \leq n
2^{\frac{-\ell +1}{2}} \sqrt{\log_2 n^2}$ indices $j$ in $\Capp (\ell)$ and likewise
for $k$, while there are at most $\frac{1}{\eps_{\ell}} +1 \leq \sqrt{\log_2 n^2} +1$
indices $i$. Hence $\# \Capp (\ell) \leq 2 n^2 2^{-\ell} (\sqrt{\log_2 n^2}+1)^3$
and (i) follows.

As for (ii), by the assumption on $R^2$ there exists $\ell \in \{0,\ldots,\lceil \log_2
\frac{n^2}{8} \rceil \}$ such that $2^{\ell -1} < R^2 \leq 2^{\ell}$.
We can now find a $B_{r_i}(j,k) \in \Capp (\ell)$ with the desired property:
Let $i$ be the largest integer such that $r_i \leq R^2$. Then by the construction of $r_i$ we have
$r_i^2/R^2 \geq 2^{-\eps_{\ell}} \geq 1-\eps_{\ell}$. Let $j$ and $k$ be the elements in 
$\{m\, d_{\ell}, m\in \N \} \cap [r_i,n-r_i+1]$ that are closest to $s$ and $t$, respectively.
Then $|j-s| \leq d_{\ell}$, $|k-t| \leq d_{\ell}$, and therefore the Euclidean distance between
$(j,k)$ and $(s,t)$ is not larger than $\sqrt{2} d_{\ell}$. Thus it follows from Lemma~\ref{CL}
below that 
\begin{align*}
|B_R(s,t) \triangle B_{r_i}(j,k)| & \leq \Bl(1-\frac{r_i^2}{R^2} +2 \frac{\sqrt{2}d_{\ell}}{R}\Br)
|B_R(s,t)|\\
& \leq \Bl(\eps_{\ell} +3 \frac{\eps_{\ell} 2^{\frac{\ell -1}{2}}}{2^{\frac{\ell -1}{2}}} \Br)
|B_R(s,t)|\\
& \leq 3 \eps_{\ell} |B_R(s,t)|\ \leq \ \frac{|B_R(s,t)|}{\sqrt{\log_2 \frac{n^2}{|B_R(s,t)|}}}.
\end{align*}
\end{proof}

\begin{lemma} \label{CL}
Let $0<r \leq R$ and $d \in \mathbb{R}^2$. Then 
$$
| B_R(0) \triangle B_r(d) | \ \leq \ \Bl(1-\frac{r^2}{R^2} +2 \frac{|d|}{R} \Br) |B_R(0)|.
$$
\end{lemma}

\begin{proof}[Proof of Lemma \ref{CL}]
\begin{align}
| B_R(0) \triangle B_r(d) | & = | B_R(0)| -| B_r(d)| +2|B_r(d) \setminus B_R(0)| \nn \\
& \leq | B_R(0)| -| B_r(d)|+2|B_R(d) \setminus B_R(0)| \nn \\
& = 3 | B_R(0)| -| B_r(d)|-2 |B_R(d) \cap B_R(0)|.  \label{L*}
\end{align}
If $|d| \leq 2R$, then $B_R(d) \cap B_R(0)$ is the union of two circular segments
with equal area. The formula for a circular segment gives

\begin{align*}
|B_R(d) \cap B_R(0)| & = 2R^2 \cos^{-1} \Bl(\frac{|d|}{2R}\Br)-|d|\sqrt{R^2 -\Bl(\frac{|d|}{2}\Br)^2}\\
& \geq 2R^2 \Bl( \frac{\pi}{2} -\frac{\pi}{2} \frac{|d|}{2R}\Br) -|d|R\\
& \geq \pi (R^2 -|d|R).
\end {align*}
Hence (\ref{L*}) is not larger than $R^2 \pi -r^2 \pi +2\pi |d|R$. The lemma follows as it 
trivially also holds in the case $|d|>2R$. 
\end{proof}

\begin{proof}[Proof of Theorem~\ref{CD}] The claims about the null distribution follow as in the
case of univariate intervals (Theorem~\ref{thm:HC-BJ-Under-the-null}) and multivariate 
rectangles (Theorem~\ref{thm:multivariate-null}). The key argument is again to show that the
balls in $\Capp (\ell)$ can be grouped into a small number of groups each consisting of
$\sim \frac{n^2}{2^{\ell+2}}$ disjoint balls. To this end, define $L_{\ell}$ to be the
smallest multiple of $d_{\ell}$ not smaller than $\max_i 2r_i$, so $L_{\ell} \sim 2\sqrt{2^{\ell}}$.
Define
$$
\mbox{shift}_{\ell} (j,k,r_i)\ :=\ \Bl\{B_{r_i}(s,t):\ s=j+uL_{\ell},\ t=k+vL_{\ell},\ 
s,t \in [r_i,n-r_i+1]; \ u,v\in \N_0 \Br\}.
$$
By construction, the balls in $\mbox{shift}_{\ell} (j,k,r_i)$ are mutually disjoint. One readily
checks
$$
\Capp (\ell)\ =\ \bigcup_{j,k \in \{d_{\ell},2d_{\ell},\ldots, L_{\ell}\},
i \in \{0,\ldots,\lfloor \frac{1}{\eps_{\ell}} \rfloor \}} \mbox{shift}_{\ell} (j,k,r_i).
$$
There are $\sim \Bl( \frac{n}{L_{\ell}}\Br)^2 \sim \frac{n^2}{2^{\ell +2}}$ balls in
$\mbox{shift}_{\ell} (j,k,r_i)$, and the number of groups is $\sim \Bl( \frac{L_{\ell}}{d_{\ell}}\Br)^2
\frac{1}{\eps_{\ell}} \sim 8 \eps_{\ell}^{-3} \leq 8 (\log n)^{\frac{3}{2}}$.
The latter number has an additional factor $(\log n)^{\frac{1}{2}}$ compared to the case
of univariate intervals, which likewise affects the convergence rate of $sHC_{n}^{(2)}$
as is clear from the proof of Theorem~\ref{thm:HC-BJ-Under-the-null}.
The proof of the optimality properties follows that of Theorem~\ref{optimality2}. 
\end{proof}

\section*{Appendix} \label{appn}

\begin{proof}[Proof of Lemma \ref{lem: goodness_of_fit_converge_to_one}]
Note that using the same considerations as in (\ref{represent}) we obtain 
\[
\sup_{n^{-4r}<p_{(i)}<n^{-4r_{0}}}\left|\frac{i}{n_{\ell^{*}}p_{(i)}}-1\right|\ =\ 
\sup_{n^{-4r}<\bar{\Phi}(t)<n^{-4r_{0}}}\left|\frac{\sum_{I\in\Iapp(\ell^{*})}
\boldsymbol{1}(\X(I)\geq t)}{n_{\ell^{*}}\bar{\Phi}(t)}-1\right|.
\]
By Lemma~\ref{lem: Short-range-dependent} the intervals
in $\Iapp(\ell^{*})$ can be grouped into $i_{\mathrm{max}}\leq 144 \log_2 n$ groups,
each of which consists of the same (up to $\pm 1$) number $N_{\ell^*}=L_n n^{1-\alp}$
of disjoint intervals as the first group.
Let $I_1,\ldots,I_{N_{\ell^*}}$ denote the intervals in the first group. Then for
$\eps \in (0,1)$
\bean
& & \P_{H_1} \left( \sup_{n^{-4r}<\bar{\Phi}(t)<n^{-4r_0}} \left| \frac{\sum_{I \in \Iapp(\ell^*)}
  {\bf 1}(\X(I) \geq t)}{n_{\ell^*} \bar{\Phi}(t)} -1 \right| > \eps \right) \nn \\
& \leq & n\ i_{\mathrm{max}} \sup_{n^{-4r}<\bar{\Phi}(t)<n^{-4r_0}} \P_{H_1} \left(
  \left| \frac{\sum_{i=1}^{N_{\ell^*}} {\bf 1}(\X(I) \geq t)}{N_{\ell^*} \bar{\Phi}(t)} 
  -1 \right| > \eps \right) \label{app*}
\eean
since there are not more than $n$ p-values in $(n^{-4r},n^{-4r_0})$. The $I_i$ being disjoint
implies that the $\X(I_i)$ are independent and that at most $2n^{1-\alp -\beta}$ of the $I_i$
can intersect with one of the $n^{1-\alp -\beta}$ intervals that have an elevated mean. Such
an overlap results in $\E \X(I) \leq \sqrt{3r \log n}$. Thus under $H_1$
\be  \label{app**}
{\rm Bin}(N_{\ell^*},\bar{\Phi}(t))\ \stackrel{d}{\leq}\ \sum_{i=1}^{N_{\ell^*}} {\bf 1}(\X(I) \geq t)\
\stackrel{d}{\leq} \ {\rm Bin}(2n^{1-\alp -\beta},\bar{\Phi}(t-\sqrt{2r \log n})) +
{\rm Bin}(N_{\ell^*},\bar{\Phi}(t))
\ee
Note that the function $\frac{\bar{\Phi}(t)}{\bar{\Phi}(t-\sqrt{2r \log n})}$ is decreasing in $t$ as can 
be seen by differentiating and employing the increasing hazard rate property of the normal distribution.
So if we define $t^*$ via $\bar{\Phi}(t^*)=n^{-4r}$, then $t^* = (1+o(1)) \sqrt{2(4r) \log n}$ and
\bea
\inf_{n^{-4r}<\bar{\Phi}(t)<n^{-4r_0}} \frac{N_{\ell^*} \bar{\Phi}(t)}{2n^{1-\alp -\beta}
  \bar{\Phi}(t-\sqrt{2r\log n})}
& = & L_n \frac{ \bar{\Phi}(t^*)}{n^{-\beta} \bar{\Phi}(t^*-\sqrt{2r\log n})}\\
& = & L_n \frac{n^{-4r}}{n^{-\beta} n^{-(\sqrt{4r}-\sqrt{r})^2}}\\
& = & L_n n^{\beta -3r}\ \rightarrow \ \infty \ \ \mbox{ as } r<\beta /3.
\eea
Hence for $n \geq n_0$ the above inf is larger than $4/\eps$ and so together with (\ref{app**})
we get
\begin{align}
\sup_{n^{-4r}<\bar{\Phi}(t)<n^{-4r_0}} {\P}_{H_1} & \left(
\frac{\sum_{i=1}^{N_{\ell^*}} {\bf 1}(\X(I) \geq t)}{N_{\ell^*} \bar{\Phi}(t)}
    > 1+\eps \right)\nn \\
& \leq  \sup_{n^{-4r}<\bar{\Phi}(t)<n^{-4r_0}} {\P} \left(
 \frac{{\rm Bin}(2n^{1-\alp -\beta},\bar{\Phi}(t-\sqrt{2r \log n}))}{N_{\ell^*} \bar{\Phi}(t)}
   > \eps/2 \right) \nn \\ 
& \hspace{0.5cm} +\sup_{n^{-4r}<\bar{\Phi}(t)<n^{-4r_0}} {\P} \left(
 \frac{{\rm Bin}(N_{\ell^*},\bar{\Phi}(t))}{N_{\ell^*} \bar{\Phi}(t)}
   > 1+\eps/2 \right)\nn \\
& \leq  \sup_{n^{-4r}<\bar{\Phi}(t)<n^{-4r_0}} {\P} \left(
 \frac{{\rm Bin}(2n^{1-\alp -\beta},\bar{\Phi}(t-\sqrt{2r \log n}))}{
2n^{1-\alp -\beta}\bar{\Phi}(t-\sqrt{2r \log n})}
   > 2 \right) \nn \\ 
& \hspace{0.5cm} +\sup_{n^{-4r}<\bar{\Phi}(t)<n^{-4r_0}} {\P} \left(
 \frac{{\rm Bin}(N_{\ell^*},\bar{\Phi}(t))}{N_{\ell^*} \bar{\Phi}(t)}
   > 1+\eps/2 \right)  \label{app***}
\end{align}
Now we use Bennett's inequality, which gives
$$
\P \left(\left| \frac{{\rm Bin}(m,p)}{mp} -1\right| > \eps \right)\ \leq \ 
2 \,\exp \left(-mp \eps^2/3 \right).
$$
Thus (\ref{app***}) is not larger than
\begin{align*}
\sup_{n^{-4r}<\bar{\Phi}(t)<n^{-4r_0}} &  2\,\exp \left( -2n^{1-\alp -\beta}
\bar{\Phi}(t-\sqrt{2r \log n})/3 \right)\\
& \hspace{0.5cm} +
\sup_{n^{-4r}<\bar{\Phi}(t)<n^{-4r_0}} 2\,\exp \left(-N_{\ell^*} \bar{\Phi}(t) \eps^2/12 \right) \\
& =  2\,\exp \left(-L_n n^{1-\alp -\beta} n^{-(\sqrt{4r} -\sqrt{r})^2}\right) +
  2\,\exp \left(-L_n n^{1-\alp} n^{-4r} \right)\\
& \leq  4 \,\exp(-L_n n^{\kappa})
\end{align*}
for some $\kappa >0$ as $r<\beta /3$ requires $\beta <\frac{3}{4}(1-\alp)$ and hence
$r<(1-\alp)/4$. The left tail probability in (\ref{app*}) is easily bounded analogously using
the left inequality in (\ref{app**}). Hence (\ref{app*}) is not larger than
$$
8\,n (144 \log_2 n)  \exp(-L_n n^{\kappa})\ \rightarrow 0.
$$
\end{proof}

\section*{Acknowledgements}

The authors were supported by NSF grants DMS-1220311 and DMS-1501767.

\bibliographystyle{apalike}
\bibliography{reference3}

\end{document}